\documentclass[]{imsart}

\RequirePackage{amsthm,amsmath}
\RequirePackage[numbers]{natbib}
\RequirePackage[colorlinks,citecolor=blue,urlcolor=blue]{hyperref}
\RequirePackage{graphicx}

\startlocaldefs
\numberwithin{equation}{section}
\theoremstyle{plain}

\usepackage{amssymb, amsfonts}
\usepackage{eucal}
\usepackage[inline,shortlabels]{enumitem}
\usepackage{graphicx}
\graphicspath{{figures/}}

\newtheorem{theorem}{Theorem}
\newtheorem{lemma}[theorem]{Lemma}

\newtheorem{corollary}[theorem]{Corollary}
\newtheorem{ass}[theorem]{Assumption}
\newtheorem{remark}[theorem]{Remark}

\newcommand{\bs}[1]{\boldsymbol{#1}}

\newcommand{\ip}[1]{\langle #1 \rangle}

\newcommand{\calN}{\mathcal N}
\newcommand{\calQ}{\mathcal Q}
\newcommand{\calS}{\mathcal S}
\newcommand{\calX}{\mathcal X}
\newcommand{\bb}[1]{\mathbb #1}
\DeclareMathOperator{\cov}{cov}
\DeclareMathOperator{\var}{var}
\DeclareMathOperator{\diag}{diag}
\DeclareMathOperator{\trace}{tr}

\newcommand{\ff}{{\bs{f\!f}}}
\newcommand{\fu}{{\bs{fu}}}
\newcommand{\uf}{{\bs{uf}}}
\newcommand{\uu}{{\bs{uu}}}
\newcommand{\fmv}{\hat f_m}
\newcommand{\kmv}{\hat k_m}

\newcommand{\tp}{\Pi(\,\cdot\,|\bs x,\bs y)}
\newcommand{\vp}{\Psi(\,\cdot\,|\bs x,\bs y)}
\newcommand{\bm}[1]{\begin{bmatrix} #1 \end{bmatrix}}
\let\P\undefined
\newcommand{\P}{\mathrm P}

\newcommand{\KL}{D_{\mathrm{KL}}}
\newcommand{\id}{\mathrm{id}}
\newcommand{\eps}{\epsilon}

\usepackage[svgnames]{xcolor}
\definecolor{ye}{cmyk}{0, .1, .8, .2}

\endlocaldefs

\begin{document}
\begin{frontmatter}
\title{Uncertainty quantification for sparse spectral variational approximations in Gaussian process regression}
\runtitle{Uncertainty Quantification for Variational Gaussian processes}

\begin{aug}
\author{\fnms{Dennis} \snm{Nieman}\ead[label=e1]{d.nieman@vu.nl}}
\address{Department of Mathematics, Vrije Universiteit Amsterdam, The Netherlands\\
\printead{e1}}

\author{\fnms{Botond} \snm{Szabo}\thanksref{t1}\ead[label=e2]{botond.szabo@unibocconi.it}}
\address{Department of Decision Sciences and BIDSA, Bocconi University, Italy\\
\printead{e2}}
\and
\author{\fnms{Harry} \snm{van Zanten}\ead[label=e3]{j.h.van.zanten@vu.nl}}
\address{Department of Mathematics, Vrije Universiteit Amsterdam, The Netherlands\\
\printead{e3}}

\thankstext{t1}{Co-funded by the European Union (ERC, BigBayesUQ, project number: 101041064). Views and opinions expressed are however those of the author(s) only and do not necessarily reflect those of the European Union or the European Research Council. Neither the European Union nor the granting authority can be held responsible for them.}

\runauthor{Nieman et al.}

\end{aug}

\begin{abstract}
We investigate the frequentist guarantees of the variational sparse Gaussian process regression model. In the theoretical analysis, we focus on the variational approach with spectral features as inducing variables. We derive guarantees and limitations for the frequentist coverage of the resulting variational credible sets. We also derive sufficient and necessary lower bounds for the number of inducing variables required to achieve minimax posterior contraction rates. The implications of these results are demonstrated for different choices of priors. In a numerical analysis we consider a wider range of inducing variable methods and observe similar phenomena beyond the scope of our theoretical findings.
\end{abstract}


\end{frontmatter}

\maketitle

\section{Introduction}
One of the key challenges in Bayesian statistics is to approximate intractable or computationally infeasible posterior distributions. This problem is becoming even more pronounced in applications where the amount of available information is rapidly growing, further increasing the complexity of the posterior. Variational methods provide a convenient way to overcome such computational issues in Bayesian statistics. The variational approximation starts by selecting an appropriate class of distributions, referred to as the variational class. Then the complex posterior distribution is approximated by its projection on this variational class with respect to the Kullback-Leibler divergence. The challenge in choosing the variational class is twofold: firstly, the variational posterior should reduce computational complexity and if possible increase the interpretability of the distribution; secondly, for meaningful inference it is crucial that the variational posterior has good statistical properties. For an overview of variational Bayes methods we refer to the review article \cite{blei2017}.

Although variational Bayes approximations are routinely used in practice, up to recently they were considered black box procedures with very limited theoretical underpinning. In the last few years the asymptotic properties of the variational posterior have been investigated. Abstract results were derived for posterior contraction rates and applied to various high-dimensional and nonparametric models, considering typically mean-field variational classes; see for instance \cite{zhang2020,yang:2020,alquier:2020,ray:2022,ray:2020}. However, almost all of these results focus on the recovery of the underlying true parameter of interest and do not address the quality of uncertainty quantification. 

In fact, one of the main appeals and strengths of the Bayesian paradigm is that it provides a probabilistic solution to the statistical problem. 
The posterior distribution can be used to quantify the remaining uncertainty about the parameter of interest. 
In practice this uncertainty is usually visualised by plotting credible regions. 
These are subsets of the parameter space with prescribed posterior probability (typically $95\%$).  In parametric models the celebrated Bernstein-von Mises theorem \cite{le2012asymptotic,van2000asymptotic} provides asymptotic frequentist coverage guarantees for credible sets under mild assumptions, meaning that 
credible sets can be interpreted as frequentist confidence sets.
In high-dimensional and non-parametric frameworks such a strong guarantee does not automatically hold in general (see e.g. \cite{cox:1993}). Nevertheless, by now we have a relatively good understanding of how to tune the prior to achieve asymptotic confidence guarantees, see for instance \cite{knapik2011,castillo:2014,szabo:2015,rousseau:2020}.

Despite its importance, so far hardly any results are available on the frequentist reliability of the variational Bayesian uncertainty quantification, where the real  credible sets are replaced by approximate credible sets derived from the variational posterior.
In fact, many of the available results are rather negative, showing that (mean-field) variational methods are often over-confident in the sense that they substantially underestimate the uncertainty of the procedure, see for instance \cite{bishop2006pattern,blei2017} for some standard examples. There are only a few positive results available. In \cite{wang2019frequentist} a variational version of the Bernstein-von Mises theorem was derived in parametric models, while in \cite{giordano2015linear} a correction was proposed using linear response methods to recover the original posterior covariance structure. However, both of these results consider only parametric models and it is in general unclear how variational credible sets behave in high- and infinite-dimensional settings. A related result we are aware of can be found in the recent paper \cite{vakili2022improved}, 
where confidence statements for a variational posterior are obtained in a special bandit-like regression setting similar to the 
results of \cite{seeger} for the true posterior. Although interesting, these results and the underlying techniques do not transfer 
to the usual nonparametric regression setting we consider in this paper.

In our analysis we focus on the popular and routinely used Gaussian process (GP) regression model. Exact computation of the posterior quickly becomes infeasible in practice since the computational cost scales cubically in the sample size. To overcome this problem, a sparse variational approximation method was proposed in \cite{titsias2009}. The variational class is parametrized by so called inducing variables which are fitted to the posterior. This approach has become increasingly popular in the machine learning community and has been applied in various settings including deep Gaussian processes and solving inverse problems. Recently theoretical guarantees were also derived for it. In \cite{burt2020} the average Kullback-Leibler distance between the variational and true posterior was studied, while in \cite{vbgp} minimax posterior contraction rates were derived. However, none of these results provided guarantees for the frequentist validity of the resulting uncertainty quantification, which is arguably one of the main aims in the Bayesian analysis.

We focus in our theoretical analysis on a specific choice of inducing variables which we call population spectral features, allowing (relatively) tractable mathematical analysis. In our numerical analysis we also show that with other choices of inducing variables similar behaviour is obtained. We observe that, in contrast to the simple parametric examples using mean-field variational approximations, in the nonparametric GP regression framework the variational posterior provides, from a frequentist perspective, reliable uncertainty statements for appropriately tuned priors. In fact, the good coverage property does not depend on the number of inducing variables used in the procedure. Besides the coverage of variational credible sets we also derive lower bounds for the number of inducing variables one has to use to achieve minimax contraction rates. This  complements the contraction rate guarantees given in \cite{vbgp}, where the sharpness of the lower bound was conjectured, but not verified. To achieve this we use a different proof technique than in \cite{vbgp}. We apply a variational version of the standard kernel ridge regression method \cite{scholkopf2002learning}. This direct approach provides more control of  posterior properties. Finally, we apply our abstract results to two GP priors with polynomially and exponentially decaying eigenvalues, respectively.

\paragraph{Contributions.} We summarize our contributions below:
\begin{itemize}
\item We give an explicit formula for the contraction rate of the variational posterior in terms of the prior and the true 
regression function. This gives a condition on the prior and the minimal number of inducing variables in the variational approximation needed to obtain minimax contraction rates.
\item If the number of inducing variables is too low, the contraction rate is sub-optimal, regardless of the choice of prior.
\item Irrespective of the number of inducing variables, variational credible sets cover the truths that are at least as smooth as the prior, whereas coverage may be bad if the prior over-smoothes the true regression function.
\end{itemize}

\paragraph{Outline.} In the next section, we describe the Gaussian process regression model studied in this paper. We recall the details of the variational procedure with inducing variables and derive a connection with kernel ridge regression, used in the proofs. Lastly we introduce the specific choice of inducing variables considered in the theoretical analysis. In Section~\ref{s:cr}, we develop a theory for contraction rates. Section~\ref{s:uq} consists of the theory on uncertainty quantification. In Section~\ref{s:examples} these results are applied to two specific priors, with polynomially and exponentially decaying eigenvalues, respectively. We conclude with a numerical analysis, including various inducing variable methods in Section \ref{sec:numerical}. The proofs are deferred to the Appendix. In Section \ref{sec:proof:thm} we prove the more abstract results, while the proofs for the examples are given in Section \ref{sec:proof:cor}.

\paragraph{Notation.} For sequences $a_n,b_n$ of non-negative real numbers, we write $a_n \lesssim b_n$ if there exists a constant $C>0$ such that $a_n \leq C b_n$ for all $n \in \bb N$. We write $a_n \asymp b_n$ if both $a_n \lesssim b_n$ and $b_n \lesssim a_n$ hold. We indicate with an apostrophe the transpose $A'$ of a matrix $A$.

\section{Variational approximations for Gaussian process regression}

Throughout the paper we investigate the nonparametric regression model
\begin{equation}\label{e:model}
y_i = f(x_i) + \varepsilon_i, \qquad i=1, \ldots, n,
\end{equation}
with i.i.d. design points $x_i$ with respect to some common probability measure $\mu$ on some $\calX \subseteq \bb R^d$, and i.i.d. mean-zero Gaussian measurement errors $\varepsilon_i \sim \calN(0,\sigma^2)$, for some $\sigma^2>0$. We view the unknown parameter $f$ as an element of the function space $L^2(\calX,\mu)$. We endow $f$ with a centered Gaussian process (GP) prior $\Pi$ determined by its covariance kernel $k: \calX \times \calX \to \bb R$.

The GP prior is conjugate for the regression model \eqref{e:model}, which means that the posterior is also a GP. However, the computational and memory costs of obtaining the posterior are $O(n^3)$ and $O(n^2)$, respectively, which is prohibitive for large data sets. Therefore, in practice various approximation methods are applied for inference, see \cite{Rasmussen2006} for a detailed discussion. One of the most commonly used approaches is the sparse variational approximation using inducing variables, proposed by \cite{titsias2009}. In the next subsection we give a brief summary of inducing variable variational Bayes methods in general and we will focus on a specific, analytically convenient version of the method, using population spectral features.

\subsection{Variational Bayes with inducing variables}
In the variational framework, the posterior distribution is approximated by projecting it onto an appropriately selected class $\calQ$ of probability measures on $L^2(\calX,\mu)$ with respect to the Kullback-Leibler divergence. Letting $\tp$ denote the true posterior, the variational posterior $\vp$ is defined as
\begin{equation}\label{e:klmin}
\vp = \arg \min_{\Psi \in \calQ} \KL(\Psi \|\tp),
\end{equation}
where $\KL$ denotes Kullback-Leibler divergence.

In the inducing variables framework, the variational class $\mathcal Q$ is constructed using a collection $\bs u = (u_1,\ldots,u_m)$ of (known, specified) bounded linear functionals evaluated at $f$. The linearity guarantees that the distribution of $\bs u$ under $\Pi$ is $m$-dimensional multivariate Gaussian. Furthermore, the prior conditional on $\bs u$ is another GP law. To obtain a low-dimensional optimization problem and preserve aspects of the prior distribution, \cite{titsias2009} proposed to fit a variational distribution of $\bs u$ to its posterior, while keeping the conditional prior distribution for $f|\bs u$. More concretely, the sparse inducing variable variational class consists of the distributions
\begin{align}
 \Psi = \int \Pi(\,\cdot \,| \bs u) \, d\Psi_{\bs u}(\bs u), \label{def: vb:class}
\end{align}
where $\Psi_{\bs u}$ is any non-degenerate $m$-dimensional Gaussian distribution, indexing the variational class $\calQ$. This way the posterior information is compressed into the fitted $m$-dimensional distribution of $\bs u$. Nevertheless, $\Psi$ is still a nonparametric distribution on $L^2(\calX,\mu)$ which is equivalent to the prior.

The variational posterior is computed by finding the distribution $\Psi_{\bs u}$ that minimizes the Kullback-Leibler divergence between $\Psi$ given in \eqref{def: vb:class} and the true posterior. As is shown in \cite{titsias2009}, a unique solution exists and can be computed analytically. The corresponding variational posterior is also GP with respective mean and covariance function
\begin{align}
\label{e:vpm}
\fmv(x) & = K_{x\bs u}(\sigma^2 K_\uu + K_\uf K_\fu)^{-1} K_\uf \bs y, \\
\label{e:vpc}
\kmv(x,y) & = k(x,y) - K_{x\bs u} K_\uu^{-1} K_{\bs u y} + K_{x\bs u}(K_\uu + \sigma^{-2} K_\uf K_\fu)^{-1} K_{\bs u y}.
\end{align}
Here $\bs y = (y_1,\ldots y_n)$ is the vector of response variables, $K_\uu$ is the covariance matrix of $\bs u$ under $\Pi$ (with entries $\cov_\Pi(u_i,u_j)$), and similarly $K_\fu = K_\uf'$ is the $n \times m$ matrix with entries $\cov_\Pi(f(x_i),u_j) = \Pi f(x_i)u_j$, and $K_{\bs u x} = K_{x \bs u}' = \cov_\Pi(\bs u, f(x))$. We note that in the special case $\bs u = \bs f$, the true posterior is recovered. The above formulas were derived in \cite{titsias2009} for the inducing points method ($u_j = f(z_j)$ for $z_j \in \calX$, $j=1,\ldots,m$), but the same computations hold for any choice of the inducing variables. For completeness we provide the details in Appendix~\ref{sec:vb:ind:var}.

The theoretical properties of this approach have been investigated in \cite{burt2020, vbgp} for various choices of inducing variables. The first paper deals with the accuracy of the variational approximation of the original posterior with respect to the Kullback-Leibler divergence. In the second paper, upper bounds were derived for the posterior contraction rate. We extend the latter result by using a different, kernel ridge regression technique allowing sharper control of the approximation. This analysis allows us to derive lower bounds for the contraction rate and to investigate the frequentist coverage properties of the credible sets resulting from the variational approximation.

\subsection{Kernel ridge regression}

The posterior mean, which is the maximum a posteriori in case of Gaussian processes, can equivalently be obtained as a kernel ridge regression (KRR) estimator. Let $\bb H$ be the reproducing kernel Hilbert space (RKHS) associated with the GP kernel $k$. Then (see e.g. \cite{kimeldorf1970}) the mean of the original posterior equals
\[ \arg \min_{f \in \bb H} \sum_{i=1}^n (y_i-f(x_i))^2 + \sigma^2\|f\|_{\bb H}^2. \]

It is not difficult to see that the variational posterior mean \eqref{e:vpm} can also be viewed as a KRR estimator, for an appropriate choice of 
the RKHS. Since the inducing variables $u_j$ are linear functionals of $f$, the functions 
\[ h_j : \calX \to \bb R,~ x \mapsto \cov_\Pi(f(x),u_j) = \Pi f(x)u_j \]
are elements of $\mathbb H$ (see \cite{rkhs}). Let $\mathbb H_m$ denote the linear subspace of $\mathbb H$ spanned by the functions $h_1,\ldots,h_m$. The variational posterior mean \eqref{e:vpm} is an element of $\bb H_m$. The following lemma states that the inducing variable variational posterior mean minimizes the same objective function as the mean of the true posterior, but over the subclass $\mathbb{H}_m\subset\mathbb{H}$.

\begin{lemma}\label{l:krr}
The variational posterior mean $\fmv$ given in \eqref{e:vpm} satisfies
\begin{equation}\label{e:vpmkrr}
\fmv = \arg \min_{f \in \mathbb H_m} \sum_{i=1}^n (y_i-f(x_i))^2 + \sigma^2 \|f\|_{\bb H}^2.
\end{equation}
\end{lemma}

For a similar result in context of the more specific Nystr\"om approximation method we refer to \cite{wild2021}. Although the above lemma holds for arbitrary choices of the inducing variables, in the upcoming sections we focus  specifically on the population spectral features approach. We believe that what we infer from our results holds more generally, as illustrated in our numerical analysis.

\subsection{Population spectral features}
In our theoretical analysis we focus on a choice of inducing variables that gives the variational posterior the interpretation of a spectral approximation to the true posterior.  We assume that the prior covariance kernel $k$ is continuous and $k \in L^2(\calX\times\calX, \mu \times \mu)$, so that it has a Mercer decomposition
\begin{equation}\label{e:mercer}
k(x,y) = \sum_{j=1}^\infty \lambda_j \varphi_j(x) \varphi_j(y),
\end{equation}
where $\lambda_j$ is a decreasing, summable sequence of nonnegative numbers and $(\varphi_j)_{j \in \bb N}$ is an orthonormal basis of $L^2(\calX, \mu)$. Under the current assumptions, the series \eqref{e:mercer} converges not only in the $L^2$ sense but also uniformly on compact subsets in the support of $\mu\times\mu$ (see \cite{steinwart2012}, Corollary 3.5). We also associate with $k$ an operator $T : L^2(\calX,\mu) \to L^2(\calX,\mu)$ given by
\begin{equation}\label{e:T}
Tg(x) = \int_\calX  k(x,y) g(y) \, d\mu(y),
\end{equation}
which is called the covariance operator. The identity \eqref{e:mercer} is equivalent to the decomposition $T = \sum_{j=1}^\infty \lambda_j \ip{\cdot,\varphi_j}\varphi_j$.  We assume that the set of functions $(\varphi_j)_{j \in \bb N}$ is uniformly bounded:

\begin{ass}\label{a:bound}
The functions $\varphi_j$ satisfy
\begin{equation}\label{e:eigbound}
\sup \{ |\varphi_j(x)| : j \in \mathbb N, x \in \calX \} =: C_{\varphi} < \infty.
\end{equation}
\end{ass}

~

In the theoretical part of this paper, we consider the inducing variables
\[ u_j = \ip{f,\varphi_j} = \int f(x) \varphi_j(x)\,d\mu(x), \qquad j =1, \ldots, m, \]
which we refer to as \textit{population spectral features}. We note that this approach requires the explicit knowledge of the basis functions $\varphi_j$. Let us write $\varphi_{1:m}(x)=(\varphi_1(x),\ldots,\varphi_m(x))$, and take $\Phi$ to be the $n\times m$ matrix whose $i$-th row is $\varphi_{1:m}(x_i)'$. It follows from Fubini's theorem that $\Pi u_iu_j = \lambda_j \delta_{ij}$ and $\Pi f(x)u_j = \lambda_j \varphi_j(x)$, so
in this case
\begin{align}
\nonumber
K_\uu &= \diag(\lambda_1,\ldots,\lambda_m) =: \Lambda, \\
\label{e:k}
K_{x \bs u} &= \varphi_{1:m}(x)'\Lambda, \\
\nonumber
K_\fu & = \Phi \Lambda.
\end{align}
In view of \eqref{e:vpm} and \eqref{e:vpc}, the variational posterior $\Psi(\,\cdot\,|\bs x,\bs y)$ is the law of a Gaussian process with mean and covariance function
\begin{align}
\label{e:vpmu}
\fmv(x) & = \varphi_{1:m}(x)'(\Lambda^{-1} + \sigma^{-2} \Phi'\Phi)^{-1} \Phi' \bs y/\sigma^2, \\
\nonumber
\kmv(x,y) & = k(x,y) - \varphi_{1:m}(x)'\Big(\Lambda - (\Lambda^{-1} + \sigma^{-2} \Phi'\Phi)^{-1} \Big)\varphi_{1:m}(y), \\
\label{e:vpcu}
& = \varphi_{1:m}(x)' (\Lambda^{-1} + \sigma^{-2} \Phi'\Phi)^{-1} \varphi_{1:m}(y) + \sum_{j=m+1}^\infty \lambda_j\varphi_j(x)\varphi_j(y),
\end{align}
where the last line follows from \eqref{e:mercer}. We note that in this setting $\bb H_m\subset \bb H$ is the linear span of the first $m$ basis functions. In the next three sections, we develop theory for the population spectral features variational posterior, which is fully characterised by \eqref{e:vpmu} and \eqref{e:vpcu}.

\section{Contraction rates}\label{s:cr}
First we investigate the asymptotic recovery property of the population spectral features variational posterior. An upper bound for the contraction rate was already derived in \cite{vbgp} for general inducing variables methods. However, the implicit approach based on GP concentration functions do not directly result in lower bounds for the contraction and therefore does not imply a lower bound on the number of inducing variables one has to apply to achieve minimax contraction rates. Therefore, in this article we take a more direct approach using kernel ridge regression techniques to understand the limitations of the variational approximation.

In this section we explicitly decompose the contraction rate to bias and variance terms, which in turn illuminates how optimal contraction imposes a condition on the minimal dimension $m$ of the approximation. It also shows why the recovery accuracy does not improve any further by increasing the dimension $m$ beyond this minimum.  Moreover, a converse result to the contraction rate statement is also given, which says that the variational posterior does not contract at the optimal rate if $m$ is too low.

\subsection{Convergence rate of the mean of the variational posterior}

We assume that the data are generated according to some true $f_0 \in L^2(\calX,\mu)$. Let us denote by $\P_0$ the measure on $\calX^n \times \mathbb R^n$ under which $(\bs x,\bs y)$ satisfies \eqref{e:model} with $f=f_0$. First, we consider the variational posterior mean $\fmv$ as an estimator of $f_0$. The next lemma decomposes the mean squared error of the estimator $\fmv$ into squared bias and variance terms of the estimator $\fmv$ under $\P_0$. The bias term $B_n$ consists of two parts. The first part is accounting for the estimation error in the subspace $\bb H_m$, while the second term equals the squared norm of the orthogonal projection of $f_0$ onto the orthogonal complement $\bb H_m^\perp$ of $\bb H_m$ in $L^2(\calX,\mu)$. The other term $W_n$ is the variance term. The proof of the lemma is deferred to Section \ref{sec:proof:loss}.

\begin{lemma}\label{l:loss}
Define
\begin{equation}\label{e:F}
\nu_j = \frac{n\lambda_j}{\sigma^2+n\lambda_j}.
\end{equation}
Let $m = m_n$ be such that $m^2n^{-1} \log n \to 0$ as $n\to\infty$. Then for any bounded $f_0 \in L^2(\calX,\mu)$,
\begin{equation}\label{e:convRate}
\P_0\|\fmv-f_0\|^2 \lesssim  B_n+W_n,
\end{equation}
where
\begin{equation}\label{def:bias:var}
B_n=\sum_{j=1}^m (1-\nu_j)^2 \ip{f_0,\varphi_j}^2 + \sum_{j>m} \ip{f_0,\varphi_j}^2\quad\text{and}\quad W_n=\frac 1n \sum_{j=1}^m \nu_j^2.
\end{equation}
\end{lemma}

\begin{remark}\label{remark:technicalm}
The mild technical assumption $m^2 n^{-1} \log n$ allows for a clean presentation. It can be weakened, but this requires the introduction of a technical term on the right-hand side of \eqref{e:convRate} as in \cite{hadji2022}.
\end{remark}

\subsection{Contraction rate of the variational posterior}

The contraction rate of the variational posterior is determined by the squared bias $B_n$ and variance $W_n$ of the posterior mean introduced above, as well as the term $V_n$ introduced below, which characterises the spread of the variational posterior.
\begin{theorem}\label{t:contraction}
Let $f_0 \in L^2(\calX,\mu)$ be a bounded function and $m=m_n\to\infty$ such that $m^2 n^{-1} \log n \to 0$. Then
\begin{equation}
 \P_0\Psi(\|f - f_0\| \geq M_n \epsilon_n \,| \bs x,\bs y) \to 0 \label{eq:thm:contraction}
 \end{equation}
for arbitrary $M_n \to \infty$, where
\[ \epsilon_n^2 = B_n+V_n+W_n, \]
with $B_n$ and $W_n$ as in \eqref{def:bias:var}, and
\begin{equation}\label{def:spread}
V_n=\frac 1n \sum_{j=1}^m \nu_j + \sum_{j>m}\lambda_j.
\end{equation}
\end{theorem}

The variance term $W_n$ is always dominated by the posterior variance term $V_n$, so it does not increase the rate of contraction. The proof of the theorem is given in Section \ref{sec:proof:thm:contraction}.

Below we investigate three terms in more detail. We present an alternative formulation which is more convenient to apply in our examples and also sheds light on how the rate depends on the dimension $m$ of the variational approximation.

\begin{remark}\label{r:effdim}
By considering separately the cases that $n\lambda_j \geq 1$ and $n\lambda_j < 1$, it follows that
\begin{equation}\label{e:nuasymp}
\nu_j \asymp 1 \wedge n\lambda_j, \qquad 1-\nu_j = \frac{\sigma^2}{n\lambda_j}\nu_j \asymp 1 \wedge (n\lambda_j)^{-1}.
\end{equation}
Let us introduce
\begin{equation}\label{e:fd}
J_n = \max\{j : n\lambda_j \geq 1\},
\end{equation}
denoting the elbow point of the above quantities. Then the terms in the contraction rate \eqref{eq:thm:contraction} can alternatively be written as
\begin{align}
\label{e:fdbias}
B_n & \asymp \sum_{j=1}^{m\wedge J_n} (n\lambda_j)^{-2}\ip{f_0,\varphi_j}^2 + \sum_{j>m \wedge J_n} \ip{f_0,\varphi_j}^2, \\
\label{e:fdvar}
V_n& \asymp \frac{m \wedge J_n}{n} + \sum_{j > m \wedge J_n} \lambda_j,\qquad W_n\asymp\frac{m \wedge J_n}{n}+ n\sum_{j = J_n+1}^{m} \lambda_j^2 .
\end{align}
These identities follow immediately from \eqref{e:nuasymp}. They show that the contraction rate does not improve any further if $m$ is increased beyond $J_n$.
\end{remark}

\subsection{Lower bounds}

Next we derive lower bounds for the variational posterior contraction rate. This in turn implies a lower bound on the number of inducing variables needed to achieve minimax recovery of the truth. Let us assume that $f_0$ belongs to the $\beta$-Sobolev space
\begin{equation}
\calS^\beta := \{f \in L^2(\calX,\mu) : \|f\|_\beta < \infty\}, \qquad \|f\|_\beta = \Big(\sum_{j=1}^\infty j^{2\beta/d} \ip{f,\varphi_j}^2\Big)^{1/2} \label{def:Sob}
\end{equation}
for some $\beta>0$. The minimax convergence rate for $f_0 \in \calS^\beta$ is $n^{-\beta/(d+2\beta)}$.

The abstract results of Theorem~\ref{t:contraction} imply that for appropriately chosen eigenvalues $\lambda_j$ the variational posterior contracts around the truth with the minimax optimal rate; see Section \ref{s:examples} for two specific examples. In both of these examples the number of inducing variables has to be at least of order $n^{d/(d+2\beta)}$, which in fact corresponds to the $L^2$-entropy of the $\beta$-Sobolev ball, also referred to as the ``effective dimension'' of the model. The same minimal dimension was obtained for various choices of priors and inducing variable methods in \cite{vbgp}. So far, there was no theoretical underpinning available for the sharpness of this threshold. The following theorem aims to fill this gap for the population spectral features variational approach: it shows optimal posterior contraction can not be achieved when $m$ is below this threshold. Concretely, if $m$ grows as a power of $n$ that is strictly smaller than $d/(d+2\beta)$, the convergence and contraction rate are strictly slower than the optimal rate $n^{-\beta/(d+2\beta)}$. The proof of the theorem is given in Section~\ref{sec:proof:badcontr}.

\begin{theorem}\label{t:badcontr}
Let $m \asymp n^r$, where $0 < r < \frac{d}{d+2\beta}$. Then there exists $f_0 \in \calS^\beta$ and $0 < p < 2\beta/(d+2\beta)$ such that
\begin{equation}\label{e:badconvergence}
\P_0 \|\fmv - f_0\|^2 \gtrsim n^{-p},
\end{equation}
irrespective of the choice of prior. Moreover, possibly along a subsequence, we have
\begin{equation}
\label{e:badcontraction}
\P_0 \Psi(\|f - f_0\|^2 \leq M n^{-p} \,| \bs x,\bs y) \to 0
\end{equation}
for any $M>0$.
\end{theorem}

\section{Uncertainty quantification}\label{s:uq}
In this section we present our main results on the frequentist validity of Bayesian uncertainty quantification resulting in from the variational approximation. To this end, let us fix $\gamma \in (0,1)$ and consider the ball 
\begin{align}
C_\gamma := \{ f : \|f-\fmv\| \leq \rho_n \},\label{def:credible}
\end{align}
 where the radius $\rho_n$ is chosen such that $\Psi(f\in C_\gamma \,| \bs x,\bs y) = 1-\gamma$. This set is referred to as the $(1-\gamma)$-credible ball of the variational posterior. In the next theorem, we first study the asymptotic size of the radius $\rho_n$ under the frequentist assumption that some $f_0 \in L^2(\calX,\mu)$ generates the data. In short, under $\P_0$, the asymptotic radius of a credible set is of the order $V_n$, which was defined in \eqref{def:spread}. This is in line with the remark made earlier that $V_n$ characterises the spread of the variational posterior. The proof is given in Appendix~\ref{sec:proof:radius}.


\begin{theorem}\label{t:radius}
Suppose that $m=m_n$ is such that $m^2n^{-1} \log n \to 0$ as $n\to\infty$. Then there exists a positive constant $C$ such that the credible ball $C_\gamma$ defined in \eqref{def:credible} has radius satisfying
\[ C^{-1} V_n \leq \rho_n^2 \leq C V_n, \]
with $\P_0$-probability tending to $1$ for any $f_0 \in L^2(\calX,\mu)$.
\end{theorem} 

We now consider the frequentist coverage of the credible set $C_\gamma$, i.e. we are interested in the probability
\begin{align} 
\P_0(f_0 \in C_\gamma) = \P_0(\|\fmv - f_0\| \leq \rho_n). \label{def:coverage}
\end{align}
The next result provides guarantees but also limitations for achieving good coverage. The proof is deferred to Section \ref{sec:proof:coverage}.

\begin{theorem}\label{t:coverage}
For $B_n, W_n$ and $V_n$ defined in \eqref{def:bias:var} and \eqref{def:spread}, respectively, consider the ratio
\begin{equation}\label{e:Rn}
R_n = \frac{B_n +W_n}{V_n}.
\end{equation}
Under the conditions of Theorem~\ref{t:radius} and assuming $f_0 \in L^2(\calX,\mu)$ is bounded,
\begin{enumerate}
\item if $R_n \to 0$, then $\P_0(\|\fmv-f_0\| \leq \rho_n) \to 1$;
\item if $R_n \lesssim 1$, then $\P_0(\|\fmv-f_0\| \leq M_n \rho_n) \to 1$ for any sequence $M_n \to \infty$;
\item if $R_n \to \infty$, then $\P_0(\|\fmv-f_0\| \leq M \rho_n) \to 0$ for any $M>0$.
\end{enumerate}
\end{theorem}

In view of \eqref{def:coverage} the above theorem presents the frequentist coverage properties of the variational credible ball. It is determined by the relation of the mean squared error $\P_0 \|\fmv-f_0\|^2$, studied in Lemma~\ref{l:loss}, and the radius $\rho_n$ of the credible set, investigated in Theorem~\ref{t:radius}. 
 The first two statements are in line with the intuition that good coverage follows from the credible set's radius being larger than the {loss}. In case the radius and loss are asymptotically comparable, good coverage can be achieved by (slightly) blowing up the credible set with a growing factor $M_n$. The third statement is the converse, i.e. if the loss exceeds the radius then coverage will be bad.

We note that in statements \textit{2} and \textit{3}, the ratio $R_n$ may be replaced by $B_n/V_n$, since the variance of the posterior mean $W_n$ is always bounded by the denominator $V_n$. The numerator $B_n$ represents the (order of the) squared bias of the variational posterior mean and the denominator $V_n$ corresponds to the variance of the variational posterior. So Theorem~\ref{t:coverage} then characterizes coverage by a comparison of bias and variance; the asymptotic coverage is good if variance dominates bias, and bad if the bias strictly dominates the variance. 

Below we demonstrate in examples that irrespective of the dimension $m$ of the variational approximation, variational credible sets will cover truths that are at least as smooth as the prior, and coverage may be bad if the prior oversmoothes the truth. In this sense the variational posterior has the same behaviour as the true posterior (see for example \cite{knapik2011,hadji2022}).

\section{Examples}\label{s:examples}

We apply our theoretical results from the previous sections to two different, commonly used choices of the eigenvalues $\lambda_j$ of the kernel in \eqref{e:mercer}. First we investigate polynomially decaying eigenvalues. We note that Mat\'ern covariance kernels (including the Ornstein-Uhlenbeck process) and Riemann-Liouville processes (including integrated Brownian motions) possess such covariance structure. As a second example we consider exponentially decaying eigenvalues. We note that the squared exponential Gaussian processes possess such exponentially decaying eigenvalues. In our examples we choose the eigenfunctions $(\varphi_j)_{j\in\mathbb{N}}$ to meet our (uniform) boundedness assumption in \eqref{e:eigbound}. This latter condition is not verified in general, but only for specific choices of covariance kernels (e.g. the Ornstein-Uhlenbeck process). Alternatively, it is possible to define a suitable GP prior using \eqref{e:mercer} by specifying $(\lambda_j)_{j \in \bb N}$ and taking any basis $(\varphi_j)_{j \in \bb N}$ that satisfies the boundedness assumption.

\subsection{Polynomially decaying eigenvalues}\label{sub:poly}

In this subsection we consider eigenvalues of the form $\lambda_j \asymp j^{-1-2\alpha/d}$, $j=1,\ldots,$ for some given $\alpha>0$, and eigenfunctions $(\varphi_j)_{j\in\mathbb{N}}$ satisfying \eqref{e:eigbound}. First, by applying Theorems~\ref{t:contraction} and \ref{t:badcontr} to the present setting, we derive contraction rates for the corresponding population spectral features variational Bayes approximations. The proof of the corollary is deferred to Section \ref{sec:proof:polycontr}.

\begin{corollary}\label{c:polycontr}
Consider the kernel $k$ in \eqref{e:mercer} with $\lambda_j \asymp j^{-1-2\alpha/d}$ for some $\alpha > d/2$ and $(\varphi_j)_{j\in\mathbb{N}}$ satisfying \eqref{e:eigbound}. Then, for $m \asymp n^r$ with $1/2>r \geq d/(d+2\alpha)$, the population spectral features variational posterior contracts around bounded functions $f_0 \in \calS^\beta$ at the rate $\epsilon_n = n^{-(\beta \wedge \alpha)/(d+2\alpha)}$, which is minimax optimal for $\alpha=\beta$. Furthermore, for $\alpha=\beta$ and $r<d/(d+2\alpha)$ there exists a $f_0\in \calS^\beta$ such that the contraction rate $\eps_n\gtrsim n^{-p}$, for some $p<\beta/(d+2\beta)$.
\end{corollary}

We note that for $1/2>r \geq d/(d+2\alpha)$ the variational approximation inherits the asymptotic recovery property of the true posterior, i.e. the true posterior contracts around $f_0\in \calS^\beta$ with rate $\epsilon_n = n^{-(\beta \wedge \alpha)/(d+2\alpha)}$, which is minimax optimal for $\alpha=\beta$, see e.g. \cite{van2008rates,knapik2011}. The characterisation of the contraction rate in Theorem~\ref{t:contraction} shows that this rate is achieved uniformly over Sobolev balls $\{f\in L^2(\calX,\mu) : \|f\|_\beta \leq C\}$. The upper bound $m\lesssim n^{r}$ with $r<1/2$ is purely technical, to satisfy the assumption $m^2n^{-1}\log n$ that has been used throughout this paper and which was justified in Remark~\ref{remark:technicalm}. In \cite{vbgp} the same contraction rate results were derived using a different proof technique without this upper bound on $m$. The lower bound $m\gtrsim n^{d/(d+2\alpha)}$, however, is sharp and of importance. Fewer inducing variables can result in sub-optimal posterior contraction around the truth, not matching the minimax rate.

Next we focus on the reliability of uncertainty quantification resulting in from the variational Bayes approximation. We apply our general Theorem~\ref{t:coverage} providing guarantees but also limitations for the frequentist coverage for the variational credible sets \eqref{def:credible}. Furthermore, in view of Theorem \ref{t:radius} we derive an upper bound for the size of the credible ball. The proof of the next corollary is given in Section \ref{sec:proof:polyuq}.

\begin{corollary}\label{c:polyuq}
Let $k$ be the kernel in  \eqref{e:mercer} with $\lambda_j \asymp j^{-1-2\alpha/d}$, for some $\alpha > 0$ and eigenfunctions $(\varphi_j)_{j\in\mathbb{N}}$, satisfying \eqref{e:eigbound}.  Consider $m \asymp n^r$ population spectral inducing variables for some $r \in (0,1/2)$.
\begin{enumerate}
\item If $\alpha\leq\beta$, then
\[ \P_0(\|\fmv-f_0\| \leq M_n \rho_n) \to 1 \]
for any bounded $f_0 \in \calS^\beta$ and any sequence $M_n \to \infty$.
\item If $\alpha<\beta$ and $r < d/(d+2\alpha)$, then 
\[ \P_0(\|\fmv-f_0\| \leq \rho_n) \to 1 \]
for any bounded $f_0 \in \calS^\beta$.
\item If $\alpha>\beta>0$, then there exists $f_0 \in \calS^\beta$ such that
\[ \P_0(\|\fmv-f_0\| \leq M \rho_n) \to 0 \]
for all $M>0$.
\end{enumerate}
Finally, for $\beta=\alpha$ and $m\geq n^{d/(d+2\alpha)}$, there exists $C>0$ such that for $f_0 \in \cal S_\beta$, we have $P_0(\rho_n\leq Cn^{-\beta/(d+2\beta)})\rightarrow 1$.
\end{corollary}

The first statement of the corollary says that by considering a prior which does not oversmooth the true function, the slightly inflated credible sets provide reliable uncertainty quantification from the frequentist perspective. In other words, by blowing up a credible ball \eqref{def:credible} with a multiplicative factor of $M_n\rightarrow\infty$, its frequentist coverage tends to one. We note that inflating the credible set is in principle equivalent to (slightly) undersmoothing the prior distribution.

In statement \textit{2}, we in fact show that such technical post-processing is not necessary by taking $\alpha<\beta$ and considering less than $m\leq n^r$, with $r < d/(d+2\alpha)$, inducing variables. We note that in view of Corollary \ref{c:polycontr}, such choices of $\alpha$ and $m$ result in an inflated contraction rate. In this case the variational posterior mean is sub-optimally far from the true $f_0$ (considering squared $L^2$-loss). At the same time the size of the credible ball is also of higher magnitude $V_n \asymp \sum_{j>m}\lambda_j$ which is large enough to dominate the expected loss of the variational posterior mean.

When considering smoother priors than the true function $f_0$, as in statement \textit{3} of the corollary, the frequentist coverage of the variational credible set will tend to zero, independently from the number of inducing variables applied.

Finally, as a follow up of Corollary  \ref{c:polycontr}, the size of the credible set achieves the minimax rate if the regularity of the prior matches the regularity of $f_0$.  All these results, in fact, are in line with the frequentist coverage results derived for the true posterior distribution \cite{knapik2011}.

\subsection{Exponentially decaying eigenvalues}

We also study the series prior with (rescaled) exponentially decaying eigenvalues, motivated by the squared exponential kernel $k(x,y) = \exp(-b\|x-y\|^2)$. In this subsection we consider eigenvalues
\begin{equation}\label{e:expdecay}
\lambda_{n,j} \asymp \exp(-\tau_n j^{1/d})\qquad j=1,\ldots,
\end{equation}
in \eqref{e:mercer} for some rescaling sequence $\tau_n\to0$. For fixed $\tau_n$, the decay of the eigenvalues is faster than in the previous example, giving the associated prior support consisting only of functions that are substantially smoother than the true function $f_0\in \calS^\beta$. To compensate for the rapidly diminishing eigenvalues a rescaling factor $\tau_n$ is introduced, which shrinks the trajectories of the prior draws, making them rougher. Let us consider the rescaling factor $\tau_n$ in the form
\begin{equation}\label{e:rescaling}
\tau_n := n^{-1/(d+2\alpha)} \log n,
\end{equation}
for some constant $\alpha > 0$. Taking $\alpha=\beta$, in fact, results in minimax contraction rate and reliable uncertainty quantification for the true posterior distribution, see for instance \cite{vaart2007, bhattacharya:pati,hadji:sz:2021}. The following corollary states that the contraction rate of the true posterior is inherited by its variational approximation for sufficiently high $m$. The proof is deferred to Section \ref{sec:proof:expcontr}.

\begin{corollary}\label{c:expcontr}
Let $k$ be the kernel in  \eqref{e:mercer} with eigenvalues satisfying \eqref{e:expdecay} with rescaling sequence \eqref{e:rescaling} for some $\alpha > 0$ and eigenfunctions $\varphi_j$, $j=1,\ldots$ satisfying \eqref{e:eigbound}. Then, for $m = m_n \geq (\tau_n^{-1} \log n)^d = n^{d/(d+2\alpha)}$ such that $m^2 n^{-1} \log n \to 0$, the variational posterior contracts around bounded $f_0 \in \cal H^\beta$ at the rate $\epsilon_n = n^{-(\beta\wedge\alpha)/(d+2\alpha)}$, which is minimax optimal for $\alpha=\beta$. Furthermore, for $\alpha=\beta$ and $m=n^{r}$, with $r<d/(d+2\alpha)$, there exists an $f_0 \in \cal H^\beta$, for which the posterior contracts at best at the sub-optimal rate $n^{-p}$, for some $p<\beta/(d+2\beta)$.
\end{corollary}

This corollary, on one hand shows that the variational posterior achieves minimax contraction rates for appropriately chosen rescaling factor and sufficiently many inducing variables $m\gtrsim n^{d/(d+2\alpha)}$. (The technical assumption $m^2 n^{-1} \log n=o(1)$ has been used throughout and was explained in Remark~\ref{remark:technicalm}) On the other hand it also implies that we need at least this many inducing variables, otherwise the variational posterior will provide sub-optimal recovery for the truth.

Next we discuss the validity of the variational Bayes uncertainty quantification. The result is slightly different in comparison with the polynomially decaying eigenvalues. The proof of the next corollary is given in Section \ref{sec:proof:expuq}.

\begin{corollary}\label{c:expuq}
Let $k$ be the kernel in  \eqref{e:mercer} with eigenvalues satisfying \eqref{e:expdecay} with rescaling sequence \eqref{e:rescaling} for some $\alpha > 0$ and eigenfunctions $\varphi_j$, $j=1,\ldots,$ satisfying \eqref{e:eigbound}. Furthermore, let the approximation have dimension $m = m_n \asymp n^r$ for some $r \in (0,1/2)$. Then
\begin{enumerate}
\item if $r \geq d/(d+2\alpha)$ and $\beta \geq \alpha > 0$, then for any bounded $f_0 \in \calS^\beta$ and any $M_n \to \infty$,
\[ \P_0(\|\fmv-f_0\| \leq M_n\rho_n) \to 1; \]
\item if $r < d/(d+2\alpha)$ then for any bounded $f_0 \in L^2(\calX, \mu)$
\[ \P_0(\|\fmv-f_0\| \leq \rho_n) \to 1; \]
\item if $\alpha>\beta>0$ and $m \geq n^{d/(d+2\alpha)}$, then there exists $f_0 \in \calS^\beta$ such that for all $M>0$
\[ \P_0(\|\fmv-f_0\| \leq M \rho_n) \to 0. \]
\end{enumerate}
\end{corollary}

Note that the only (real) difference between this result and Corollary \ref{c:polyuq} is that in the second statement no assumption is made about the (Sobolev) smoothness of $f_0$. If the dimension of the approximation is of order $m_n = o(n^{d/(d+2\alpha)})$, then asymptotically the coverage of the variational credible set will always be good. The rest of the results are the same as in Corollary \ref{c:polyuq}, except of the minor difference that in the third statement the lower bound on $m$ holds exactly and not up to a multiplicative constant due to the exponentially decaying form of the eigenvalues, hence the conclusions are the same as well.



\section{Numerical experiments}\label{sec:numerical}

In this section we demonstrate how the developed theory can be applied in practice. Moreover, we show using synthetic data sets that the variational Bayes method proposed by \cite{titsias2009} provides reliable uncertainty quantification (from a frequentist perspective) for appropriately tuned Gaussian process priors, independently from the number of inducing variables applied.

\subsection{Synthetic data set}

In a simulation study, we go beyond the population spectral features inducing variable method considered in our theoretical analysis above, and include other, practically more advantageous and popular approaches as well. We aspire to extrapolate our findings about the principles that govern the method with the investigated particular choice of inducing variables to other choices of inducing variables, for which the theoretical analysis is more complicated. With this goal in mind, in the current section, we aim not only to illustrate our findings but also compare various methods, to study empirically if and how our theoretical results may carry over. Between each of three choices of inducing variables, we shall point out what are common features and what is different. Besides the population spectral features, we shall also study inducing points and what we call `empirical spectral features'.

In inducing points methods, the choice of inducing variables is $u_j = f(z_j)$ for a given set of inducing points $z_1,\ldots,z_m \in \calX$. In the literature various approaches were proposed to choose the inducing points. In our numerical analysis we consider two specific choices. First we consider the equidistant choice of the inducing points, i.e. we simply take the $z_j$ on an equispaced grid on $[-\pi,\pi]$. This approach is designed to explore all parts of the signal equally. In case the design points $x_1,\ldots,x_n$ are not uniformly distributed the grid can be modified accordingly to mimic the underlying distribution of the covariates. We also consider the finite fixed-size determinantal point process ($m$-DPP), see \cite{GPBV19}.

The other choice of inducing variables is the sample analogue of the population spectral features $u_j = \ip{f,\varphi_j}$. Instead of diagonalizing the covariance operator $T$ in \eqref{e:T}, we decompose the covariance matrix of the prior at the design points $K_\ff = [k(x_i,x_j)]_{1 \leq i,j \leq n}$. Since $K_\ff$ is positive definite, there exists an orthonormal set of eigenvectors $\bs v_1,\ldots,\bs v_n$ such that $ K_\ff = \sum_{j=1}^n \mu_j \bs v_j \bs v_j'$, 
where $\mu_1 \geq \mu_2 \geq \cdots \geq \mu_n \geq 0$ are the eigenvalues of $K_\ff$. The \textit{empirical spectral features} are defined $u_j := \bs v_j' \bs f = \sum_{i=1}^n v_{j,i} f(x_i)$. In \cite{vbgp} it was shown that the empirical spectral features induce a variational posterior with similar contraction rate results as the population spectral features under similar threshold for the number of inducing variables. Furthermore, the approximation accuracy of the above variational Bayes methods to the true posterior distribution with respect to the expected Kullback-Leibler divergence was studied in \cite{burt2020}.

~

We consider the function space $L^2([-\pi,\pi],\mu)$, where $\mu$ is taken as the uniform distribution on $[-\pi,\pi]$. As our underlying true function (plotted in black in the figures) we take 
\[ f_0 = \sum_{\ell=1}^\infty \varphi_{3\ell} (3\ell)^{-1/2-\beta}/\log(3\ell) \in \mathcal H^\beta \]
with  $\beta = 0.5$, where $\varphi_j$ denote the standard Fourier basis. We use a synthetic data set with $n=2500$ realisations $x_i\sim_{\textrm{i.i.d.}} \mu$, $i=1,\ldots,n$ and independent data points $y_i\sim_{\text{ind}}\mathcal N(f_0(x_i),\sigma^2)$ with $\sigma = 0.1$. 

\begin{figure}
\includegraphics[width=5.5cm]{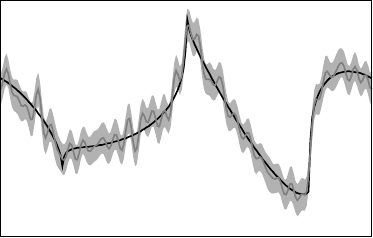}\hspace{.3cm}
\includegraphics[width=5.5cm]{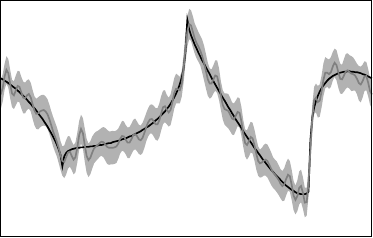}
\caption{The true posterior distribution corresponding to the GP prior with polynomially (left) and exponentially (right) decaying eigenvalues}\label{fig:posterior}
\end{figure}

We consider both the polynomially and exponentially decaying eigenvalues $\lambda_j$ from Section~\ref{s:examples} and take the GP prior with covariance kernel $k(x,y) = \sum_{j=1}^\infty \lambda_j \varphi_j(x)\varphi_j(y).$ Accordingly, in Figure \ref{fig:posterior} we plot the mean (dark gray) and 95\% pointwise credible sets of the associated true posterior when  $\lambda_j = j^{-1-2\beta}$ and $\lambda_j = \tau_n \exp(-\tau_nj/4)$, with $\tau_n = n^{-1/(1+2\beta)} \log n$, respectively.

\begin{figure}[p]
\includegraphics[width=5.5cm]{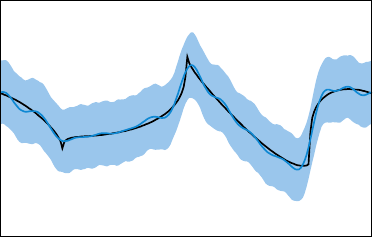}\hspace{.3cm}
\includegraphics[width=5.5cm]{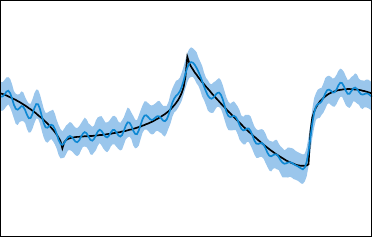}\\
\vspace{0.3cm}
\includegraphics[width=5.5cm]{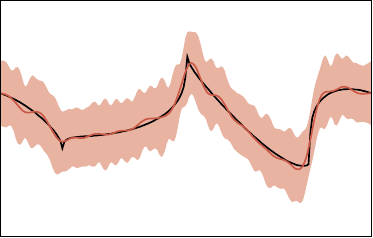}\hspace{.3cm}
\includegraphics[width=5.5cm]{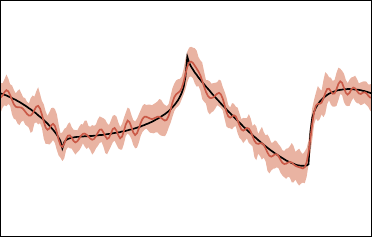}\\
\vspace{0.3cm}
\includegraphics[width=5.5cm]{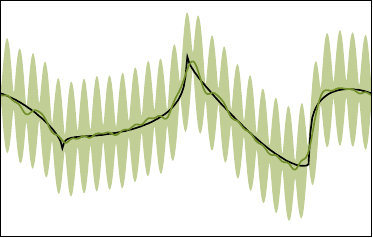}\hspace{.3cm}
\includegraphics[width=5.5cm]{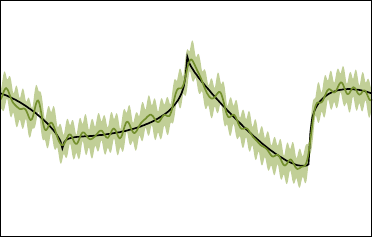}\\
\vspace{0.3cm}
\includegraphics[width=5.5cm]{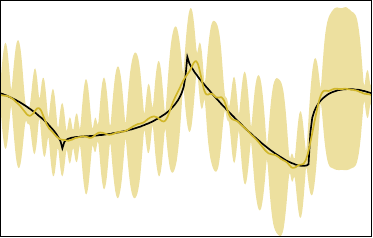}\hspace{.3cm}
\includegraphics[width=5.5cm]{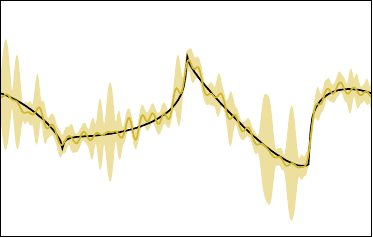}
\caption{The variational Bayes approximations for GP with polynomially decaying eigenvalues.  for $m=30$ (left) and $m=60$ (right) inducing variables. The choices of inducing variable methods from top to bottom: population spectral features (blue), empirical spectral features (red), inducing variables with equidistant design (green) and with $m$-DPP (yellow). In each figure the true function is plotted in black, the posterior mean is drawn with solid colored curves, while the shaded area illustrate the $95\%$ point-wise credible bands.}\label{fig:vb:poly}
\end{figure}

\begin{figure}[p]
\includegraphics[width=5.5cm]{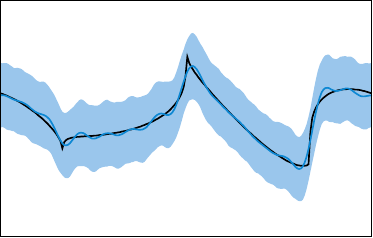}\hspace{.3cm}
\includegraphics[width=5.5cm]{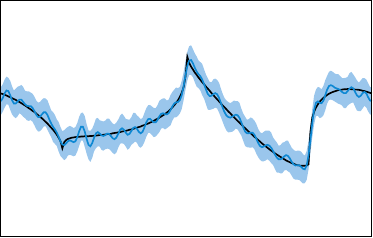}\\
\vspace{.3cm}
\includegraphics[width=5.5cm]{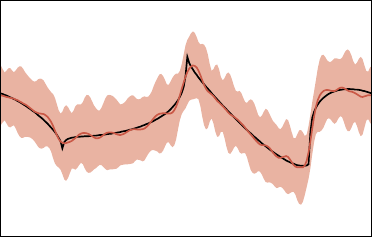}\hspace{.3cm}
\includegraphics[width=5.5cm]{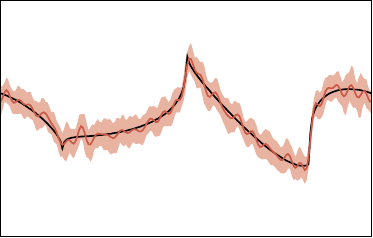}\\
\vspace{.3cm}
\includegraphics[width=5.5cm]{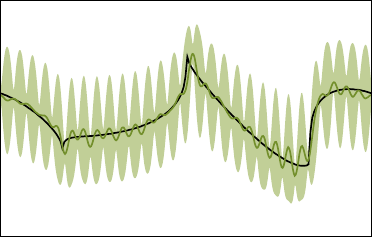}\hspace{.3cm}
\includegraphics[width=5.5cm]{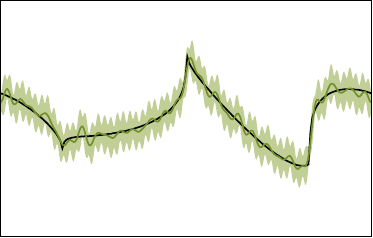}\\
\vspace{.3cm}
\includegraphics[width=5.5cm]{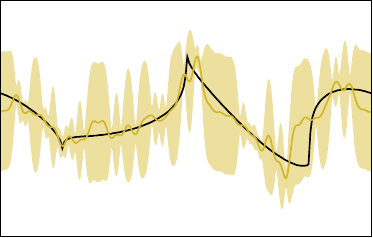}\hspace{.3cm}
\includegraphics[width=5.5cm]{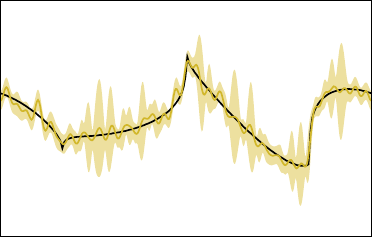}
\caption{The variational Bayes approximations for GP with exponentially decaying eigenvalues.  for $m=30$ (left) and $m=60$ (right) inducing variables. The choices of inducing variable methods from top to bottom: population spectral features (blue), empirical spectral features (red), inducing variables with equidistant design (green) and with $m$-DPP (yellow). In each figure the true function is plotted in black, the posterior mean is drawn with solid colored curves, while the shaded area illustrate the $95\%$ point-wise credible bands.}\label{fig:vb:exp}
\end{figure}

First we investigate the variational Bayes approximations for the posterior corresponding to the polynomial GP priors in Figure \ref{fig:vb:poly}. We plot the population spectral feature variational approximation in the first line (in blue), the empirical spectral feature approach in the second line (in red), the equidistant inducing points method in the third line (in green) and the $m$-DPP inducing points method in the fourth line (in yellow). The black curve stands for the underlying true function $f_0$, the colored curve for the posterior mean, while the shaded area represents the $95\%$ pointwise credible bands. We consider two choices for the number of inducing variables. For $m=30$ (left hand side) we are below, while for $m=60$ (right hand side) we are above the theoretical threshold $m=n^{1/(1+2\beta)}=50$ obtained in our analysis. One can observe that for all methods the size of the credible bands are overly large for $m=30$, while for $m=60$ they closely resemble the true posterior given on the left hand side of Figure \ref{fig:posterior} (with perhaps the exception of the $m$-DPP method). Furthermore, for both choices of $m$ the credible bands contain the true functional parameter $f_0$, illustrating the good frequentist coverage properties of the variational methods, in line with our theoretical findings. One can also notice that the approximations from the population and empirical spectral features are fairly similar. For the equidistant  inducing points method the variational posterior means do not seem to differ significantly of those in the preceding plots, but credible regions look fairly different. At the same time in the $m$-DPP method the posterior in certain neighbourhoods can be quite different from the true posterior. In case of both inducing points methods the intervals are narrow near those points on the horizontal axis that are close to the inducing points, and widen as the distance to an inducing point increases. Here, too, the credible regions seem to be over-conservative, i.e. their width being about equal to that of the credible set of the true posterior at inducing points and larger elsewhere.

Then in Figure \ref{fig:vb:exp} we plot the posterior means and credible sets resulting in from exponentially decaying eigenvalues $(\lambda_j)_{j\in\mathbb{N}}$. We use the same experimental setup as for polynomial eigenvalues, i.e. we take $m=30$ and $m=60$ inducing variables on the left and right hand side of the figure, respectively, and consider the above discussed four variational approximations. The plots are rather similar to what whas shown for the polynomial eigenvalues in Figure \ref{fig:vb:poly}, and hence the same conclusions hold for this prior.

\newpage

\subsection{Real world data}

To illustrate how the procedure can be applied in practice, we consider a real data set consisting of hourly tin oxide measurements (the PT08.S1 series from \cite{misc_air_quality_360} used to predict carbon monoxide). Since the series exhibits periodicity, a series prior with Fourier basis is suitable. With the population spectral features variational procedure we solve the time series problem of estimating the trend and periodic components.

We start by preprocessing the data. The missing observations are estimated by interpolation of neighbouring days. This introduces bias, which however we do not account for, as the main focus here is on the practical application of the considered variational GP approach. We select the first $n=9240$ observations, corresponding to $55$ weeks of data.

In the preceding synthetic examples, we used the standard ordering of the Fourier basis, meaning the periods are sorted descendingly. The same approach here would result in overly small mass on the basis functions associated with the daily and weekly periodic behaviour. Hence we reorder the basis in a data-driven way, sorting the first $n$ basis functions according to the size of $\varphi_j(\bs x)' \bs y/n$, which estimates the coefficient $\ip{\varphi_j, f_0}$.

In the Bayesian analysis we consider the series prior with covariance kernel
\[ k(x,y) = \sum_{j=1}^\infty \tau j^{-1-2\alpha} \varphi_{(j)}(x)\varphi_{(j)}(y), \]
where the subscript between brackets $(j)$ indicates the reordered indexes. We also introduced a rescaling factor $\tau>0$ for the prior eigenvalues for additional flexibility. We fix the regularity parameter to be $\alpha = 1/2$ (the `roughest' prior allowed in our theory) and the scaling parameter $\tau=2\cdot10^5$ in our experiment. These quantities in practice are typically taken in a data driven, adaptive way. However, in this work we do not address adaptation and leave its theoretical understanding for future work.

In the synthetic examples, we considered the error variance $\sigma^2$ fixed. In the real data set it needs to be estimated. A canonical way to do this is by an empirical Bayes procedure, maximising the evidence
\[ \hat\sigma = \arg\max_{\sigma} p_{\sigma}(\bs x, \bs y) = \arg\max_{\sigma}\int p_{f,\sigma}(\bs x,\bs y) \,d\Pi(f). \]
Alternatively, as proposed by \cite{titsias2009}, in the variational framework, variance estimation can be done by maximising the evidence lower bound (ELBO)
\begin{multline*}
\log \int \exp\Big(\int \log p_{f,\sigma}(\bs x,\bs y) \,d\Pi(f|\bs u)\Big) d\Pi_{\bs u}(\bs u) = \frac 12 \log |\Sigma^{-1} K_\uu| \\
 - \frac{n}{2}\log 2\pi\sigma^2 - \frac{1}{2\sigma^2} \bs y'[I - K_\fu \Sigma^{-1} K_\uf/\sigma^2]\bs y  - \frac{1}{2\sigma^2} \trace(K_\ff - K_\fu K_\uu^{-1} K_\uf),
\end{multline*}
where $\Sigma = K_\uu + \sigma^{-2} K_\uf K_\fu$ (for more details we refer to equation (28) in \cite{titsias2009} and our Appendix~\ref{sec:vb:ind:var}). This provides a significant reduction in computation time, since optimisation of the Bayes evidence requires repeated computation of the inverse of an $n\times n$ matrix, whereas the optimisation of the ELBO requires computation of the much smaller $m \times m$ inverse of $\Sigma$ (note that the matrix $K_\uu$ is $m \times m$ diagonal with entries $\tau j^{-1-2\alpha}$). 

\begin{figure}[h]
\includegraphics{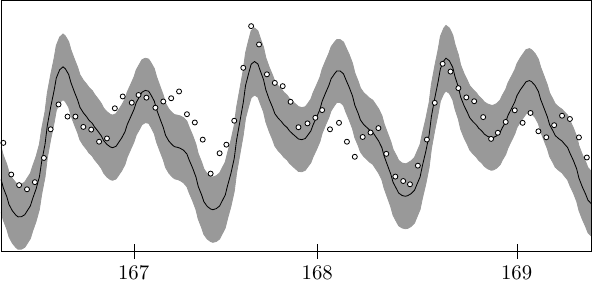}\\
~\\
\includegraphics{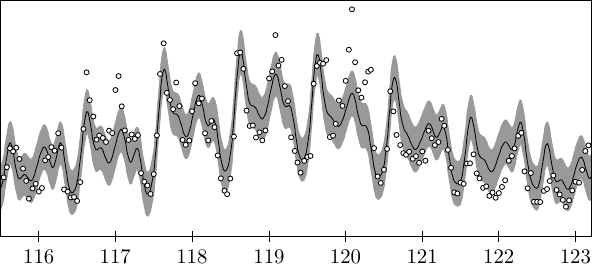}
\caption{Variational posterior on two intervals of different lengths. On the horizontal axis is the elapsed time in days. We used $m=100, \tau=2\cdot10^5$ and using the ELBO estimated $\hat\sigma \approx124$.}\label{f:realm100}
\end{figure}

Figure~\ref{f:realm100} shows the data and variational posterior mean (black line) and 95\% pointwise credible sets (gray) for our procedure with $m = 100 \approx n^{1/(1+2\alpha)}$, which is the recommendation that follows from the contraction rate results. We note that our theoretical results give frequentist coverage guarantees for the $L^2$-credible sets (which are harder to visualize) under the assumption that  $f_0\in\mathcal{S}^{\alpha}$, for $\alpha=1/2$. Nevertheless, together with the simulation study using the synthetic data set, it gives an indication of the reliability of the Bayesian uncertainty quantification from a  frequentist perspective.

\vfill

\newpage \appendix

\section{Proof of the general theorems and lemmas}\label{sec:proof:thm}

\subsection{Proof of Theorem \ref{t:contraction}}\label{sec:proof:thm:contraction}

In view of Lemma~\ref{l:loss} and the inequality
\[ \P_0 \Psi(\|f - f_0\|^2 \geq (M_n \epsilon_n)^2 \,| \bs x,\bs y) \lesssim \frac{\P_0 \Psi[\|f-\fmv \|^2 | \bs x,\bs y] + \P_0 \|\fmv - f_0\|^2}{(M_n\epsilon_n)^2} \]
(by Markov's inequality) it suffices to establish $\P_0 \Psi[\|f-\fmv \|^2 | \bs x,\bs y] \lesssim \epsilon_n^2$. Using the formula for the posterior covariance \eqref{e:vpcu}, it follows that
\begin{align}
\nonumber
\P_0 \Psi[\|f-\fmv \|^2 \,| \bs x,\bs y] 
\nonumber
& = \P_0 \int \kmv(x,x) \, d\mu(x) \\
\label{e:normSq}
& = \P_0 \trace ((\Lambda^{-1} + \sigma^{-2}\Phi'\Phi)^{-1}) + \sum_{j=m+1}^\infty \lambda_j.
\end{align}
The proof is completed by showing
\begin{equation}\label{e:learningcurve}
\P_0 \trace((\Lambda^{-1} + \sigma^{-2}\Phi'\Phi)^{-1}) \lesssim \frac 1n \sum_{j=1}^m \nu_j.
\end{equation}
To do so, the expectation term is distributed over the event
\[ \Omega_m := \Big\{ \sup_{1\leq j,k\leq m} \Big[\frac1n\sum_{i=1}^n\varphi_j(x_i)\varphi_k(x_i)-\delta_{jk}\Big]^2 > C \frac{\log n}{n} \Big\} \]
and its complement. Note that 
\[ \trace((\Lambda^{-1}+\sigma^{-2}\Phi'\Phi)^{-1}) \leq \trace \Lambda \leq \sum_{j=1}^\infty \lambda_j < \infty, \]
hence by the union bound and Hoeffding's inequality,
\[ \P_0 1_{\Omega_m} \trace((\Lambda^{-1} + \sigma^{-2}\Phi'\Phi)^{-1}) \lesssim \P_0(\Omega_m) \leq \sum_{j,k=1}^m 2 \exp\Big( \frac{-2n^2 C\log n}{4C_\varphi^2 n^2} \Big) = 2m^2 n^{-C/2C_\varphi^2}, \]
where $C$ can be chosen arbitrarily large. This term is dominated by the following upper bound that we establish below:
\begin{equation}\label{e:trbd}
\P_0 1_{\Omega_m^c} \trace((\sigma^2 \Lambda^{-1} + \Phi'\Phi)^{-1}) \lesssim \trace ((\sigma^2 \Lambda^{-1} + n I)^{-1}) = \frac{1}{n} \sum_{j=1}^m \nu_j.
\end{equation}

Using the matrix identity $A^{-1} - B^{-1} = A^{-1}(B-A)B^{-1}$ it follows that
\begin{align*}
& |\trace((\sigma^2 \Lambda^{-1} + \Phi'\Phi)^{-1}) - \trace((\sigma^2 \Lambda^{-1} + nI)^{-1})| \\
& = \trace\Big( (\sigma^2 \Lambda^{-1} + \Phi'\Phi)^{-1} (nI-\Phi'\Phi) (\sigma^2 \Lambda^{-1} + nI)^{-1}\Big) \\
& \leq \|(\sigma^2 \Lambda^{-1} + \Phi'\Phi)^{-1}\| \|nI - \Phi'\Phi\| \trace((\sigma^2 \Lambda^{-1} + nI)^{-1}),
\end{align*}
where $\|A\|$ denotes the operator norm of the matrix $A \in \mathbb R^{m \times m}$ with respect to the Euclidean norm on $\mathbb R^m$, which for symmetric $A$ coincides with the largest absolute eigenvalue of $A$. We now show that on $\Omega_m^c$, the product of the two norms in the preceding display is $o(1)$, implying \eqref{e:trbd}.

By Lemma~\ref{l:circle} below, there exists $j$ between $1$ and $m$ such that the smallest eigenvalue of $\sigma^{2}\Lambda^{-1} + \Phi'\Phi$ (on $\Omega_m^c$) is bounded from below by 
\begin{align*}
& \sigma^2 \lambda_j^{-1} + [\Phi'\Phi]_{j,j} - \sum_{\substack{1 \leq k \leq m \\ k \neq j}} |[\Phi'\Phi]_{j,k}| \\
& = \sigma^2 \lambda_j^{-1} + \Big[\sum_{i=1}^n \varphi_j(x_i)^2 - \sum_{\substack{1 \leq k \leq m \\ k \neq j}} \Big|\sum_{i=1}^n \varphi_k(x_i)\varphi_j(x_i)\Big|\Big] \\
& \geq n  - m\sqrt{C n \log n} \geq c n,
\end{align*}
where the latter inequality being true for some (deterministic) $c>0$ provided $n$ is large enough, as follows from the assumption $m^2n^{-1}\log n \to 0$. Consequently, the largest eigenvalue of the inverse $(\sigma^2\Lambda^{-1}+\Phi'\Phi)^{-1}$ satisfies
\[ \|(\sigma^2 \Lambda^{-1} + \Phi'\Phi)^{-1}\| \leq \frac{1}{cn}. \]
Another application of Lemma~\ref{l:circle} similarly shows that
\[ \|nI - \Phi'\Phi\| \leq m\sqrt{Cn \log n}. \]
Since $m \sqrt{n^{-1} \log n} \to 0$, the product of the above two norms vanishes (deterministically) on $\Omega_m^c$, implying \eqref{e:trbd} and concluding the proof.

\begin{lemma}[Gershgorin circle theorem; see \cite{gershgorin1931}]\label{l:circle}
Let $A \in \bb R^{n \times n}$ with entries $A_{ij}$. For any eigenvalue $\lambda$ of $A$, there exists $j \in \{1,\ldots,n\}$ such that
\[ |\lambda - A_{jj}| \leq \sum_{i \neq j} |A_{ij}|. \]
\end{lemma}

\subsection{Proof of Theorem \ref{t:badcontr}}\label{sec:proof:badcontr}
Take
\[ f_0 = \sum_{j=1}^\infty j^{-(1+p/r)/2} \varphi_j, \qquad\text{for some } \frac{2r\beta}{d} < p < \frac{2\beta}{d+2\beta} \]
and note that $f_0 \in \calS^\beta$ and 
\begin{equation}\label{e:bigbias}
\sum_{j>m} \ip{f_0,\varphi_j}^2 \asymp \int_m^\infty x^{-1-p/r}\,dx \asymp m^{-p/r} \asymp n^{-p} \gg n^{-2\beta/(d+2\beta)}. \end{equation}
The left-hand side in this display is a lower bound for $\|\fmv-f_0\|^2$ since $\ip{\fmv,\varphi_j} = 0$ for $j>m$, so \eqref{e:badconvergence} follows.

For the contraction rate statement we fix $p < p' < \frac{2\beta}{d+2\beta}$ and make a case distinction, first assuming
\begin{equation}\label{e:smallLambda}
\sum_{j>m} \lambda_j \lesssim n^{-p'}.
\end{equation}
Note that given $M>0$, by \eqref{e:bigbias} there exists $C_0>0$ such that for $n$ large enough
\[ Mn^{-p'} \leq (1-C_0) \sum_{j>m} \ip{f_0,\varphi_j}^2. \]
Combining this with
\[ \|f-f_0\|^2 = \sum_{j=1}^\infty \ip{f-f_0,\varphi_j}^2 \geq \sum_{j>m} \ip{f_0,\varphi_j}^2 - 2 \sum_{j>m} \ip{f_0,\varphi_j}\ip{f,\varphi_j} \]
yields
\begin{align*}
& \Psi(\|f-f_0\|^2 \leq M n^{-p'} \,| \bs x,\bs y) \\
& \leq \Psi\Big(\sum_{j>m} \ip{f_0,\varphi_j}^2 - Mn^{-p'} \leq 2 \sum_{j>m} \ip{f_0,\varphi_j}\ip{f,\varphi_j} \,\Big|\, \bs x,\bs y \Big) \\
& \leq \Psi\Big(\sum_{j>m} \ip{f_0,\varphi_j}^2 \leq 2 C_0 \sum_{j>m} \ip{f_0,\varphi_j}\ip{f,\varphi_j} \,\Big|\, \bs x,\bs y \Big) \\
& = \frac 12 \Psi\Big(\sum_{j>m} \ip{f_0,\varphi_j}^2 \leq  2C_0 \big|\sum_{j>m} \ip{f_0,\varphi_j}\ip{f,\varphi_j}\big| \,\Big|\, \bs x,\bs y \Big),
\end{align*}
where in the last line we used that under the variational posterior the inner products $\ip{f,\varphi_j}, j>m$ have a mean-zero Gaussian distribution. Continuing the preceding display, and using the Markov, triangle and Cauchy-Schwarz inequalities, it follows that
\begin{align*}
& \Psi(\|f-f_0\|^2 \leq M n^{-p'} \,| \bs x,\bs y) \\
& \lesssim \frac{\sum_{j>m} \ip{f_0,\varphi_j} \Psi(|\ip{f,\varphi_j}| \,| \bs x,\bs y)}{\sum_{j>m}\ip{f_0,\varphi_j}^2} \\
& \lesssim n^p \sum_{j>m} \ip{f_0,\varphi_j} \sqrt{\lambda_j} \\
& \leq n^p \sqrt{\sum_{j>m}\ip{f_0,\varphi_j}^2 \sum_{j>m}\lambda_j} \\
& \lesssim n^{(p-p')/2} \to 0,
\end{align*}
where in the last line we used \eqref{e:bigbias} and \eqref{e:smallLambda}.

If the assumption \eqref{e:smallLambda} does not hold, then along a subsequence,
\begin{equation}\label{e:largeLambda}
n^{p'} \sum_{j>m} \lambda_j \to \infty.
\end{equation}
In this case, we note that for $j>m$ the inner products $\ip{f,\varphi_j}\stackrel{d}{=}\sqrt{\lambda_j} Z_j$ with $Z_{j}\sim_{\text{i.i.d.}} N(0,1)$, hence by Markov's inequality,
\begin{align*}
& \Psi(\|f\|^2 \leq Mn^{-p'} \,| \bs x,\bs y) \\
& \leq \Psi\Big(\sum_{j>m} \ip{f,\varphi_j}^2 \leq M n^{-p'} \,\Big|\, \bs x,\bs y \Big) \\
& = \Psi\Big(\exp(-n^{p'}{\textstyle\sum_{j>m} \ip{f,\varphi_j}^2}) \geq \exp(-M) \,\Big|\, \bs x,\bs y \Big) \\
& \lesssim \prod_{j>m} \bb E \exp\Big( - n^{p'} \lambda_j Z_j^2 \Big) 
 = \prod_{j>m} (1+2n^{p'}\lambda_j)^{-1/2}\\
& \leq \Big(1+\sum_{j>m}2n^{p'}\lambda_j\Big)^{-1/2},
\end{align*}
which vanishes along a subsequence by the assumption \eqref{e:largeLambda}. We conclude that in this case the assertion of the theorem holds for $f_0 = 0$.

\subsection{Proof of Theorem \ref{t:radius}}\label{sec:proof:radius}
In view of Markov's inequality and the definition of $\rho_n$
\begin{equation}\label{e:rUpper}
\rho_n^2 = \rho_n^2 \gamma^{-1} \Psi(\|f - \fmv\| > \rho_n \,|\bs x,\bs y) \leq \gamma^{-1} \Psi(\|f-\fmv\|^2\,| \bs x,\bs y),
\end{equation}
Then the upper bound for $\rho_n$ is implied by the preceding display together with the inequalities \eqref{e:normSq} and \eqref{e:learningcurve} in the proof of Theorem~\ref{t:contraction}.

Next we establish the lower bound. First note, that in view of Fubini's theorem and the expansion \eqref{e:vpcu} of the variational posterior covariance,
\begin{align*}
& \Psi(\ip{f-\fmv,\varphi_i}\ip{f-\fmv,\varphi_j} \,| \bs x,\bs y) \\
& = \Psi\Big( \int (f-\fmv)\varphi_i \,d\mu \int (f-\fmv)\varphi_j \,d\mu \,\Big|\,\bs x,\bs y\Big) \\
& = \int \int \kmv(u,v) \varphi_i(u) \varphi_j(v) d\mu(u)d\mu(v) \\
& = \begin{cases}
[(\Lambda^{-1} + \sigma^{-2}\Phi'\Phi)^{-1}]_{i,j} &\text{for $i\vee j \leq m$,} \\
\lambda_j \delta_{ij} & \text{for $i\vee j > m$.}
\end{cases}
\end{align*}
Therefore, under the variational posterior $\Psi(\,\cdot\,|\bs x,\bs y)$ the inner products $\ip{f-\fmv,\varphi_j}$ are centered Gaussian random variables, where the vector of the first $m$ variables has covariance matrix $(\Lambda^{-1} + \sigma^{-2}\Phi'\Phi)^{-1}$ and is independent of the remaining $j>m$, which are mutually independent with variance $\lambda_j$. Furthermore, note that if $Z$ has an $m$-dimensional normal distribution with mean zero and covariance matrix $A$ with eigenvalues $a_1,\ldots,a_m$, then
\begin{align*}
\bb E e^{-v \|Z\|^2} &= \bb E e^{-v\sum_{j=1}^m Z_i^2}\\
& =\det(I + 2v A)^{-1/2}\\
& = \prod_{j=1}^m (1+2v a_j)^{-1/2} \\
& \leq \Big(1+2v\sum_{j=1}^m a_j\Big)^{-1/2}\\
& = \Big(1+2v \trace(A)\Big)^{-1/2} . 
\end{align*}
Therefore, by using Markov's inequality,
\begin{align*}
1-\gamma 
& = \Psi(\|f-\fmv\| \leq \rho_n \,| \bs x,\bs y) \\
& = \Psi\big(\exp(-\rho_n^{-2}\|f-\fmv\|^2) \geq \exp(-1) \,| \bs x,\bs y\big) \\
& \leq e \Psi\Big(\exp\Big(-\rho_n^{-2}\sum_{j=1}^\infty \ip{f-\fmv,\varphi_j}^2 \Big) \,\Big|\, \bs x, \bs y\Big) \\
& = e \Psi\Big(\exp\Big(-\rho_n^{-2}\sum_{j=1}^m \ip{f-\fmv,\varphi_j}^2 \Big) \,\Big|\, \bs x,\bs y\Big) \cdot \prod_{j>m} \Psi\Big(\exp(-\rho_n^{-2}\ip{f,\varphi_j}^2) \,\Big|\, \bs x,\bs y\Big) \\
& = e \det\Big(I + 2\rho_n^{-2} (\Lambda^{-1}+\sigma^{-2}\Phi'\Phi)^{-1}\Big)^{-1/2} \cdot \prod_{j>m}(1+2\rho_n^{-2}\lambda_j)^{-1/2} \\
& \leq e \Big(1 + 2\rho_n^{-2}\trace\big((\Lambda^{-1}+\sigma^{-2}\Phi'\Phi)^{-1}\big) + 2\rho_n^{-2}{\textstyle\sum_{j>m}\lambda_j}\Big)^{-1/2},
\end{align*}
which in turn implies that
\[ \rho_n^2 \geq \frac{2}{e^2(1-\gamma)^{-2}-1} \Big( \trace\big((\Lambda^{-1}+\sigma^{-2}\Phi'\Phi)^{-1}\big) + {\textstyle\sum_{j>m}\lambda_j} \Big). \]
The lower bound in the statement of the theorem now follows by the inequality (12) from \cite{vivarelli1999}, which reads
\[ \trace\big((\Lambda^{-1} + \sigma^{-2}\Phi'\Phi)^{-1}\big) \geq \sum_{j=1}^m \frac{ \sigma^2\lambda_j}{\sigma^2+ \lambda_j  \sum_{i=1}^n \varphi_j(x_i)^2} \geq (C_\varphi^{-2} \wedge 1) \sigma^2 \frac 1n \sum_{j=1}^m \nu_j. \qedhere \]

\subsection{Proof of Theorem \ref{t:coverage}}\label{sec:proof:coverage}
Consider first the statements \textit{1} and \textit{2}. Note that, with $K$ as in Theorem~\ref{t:radius}, by Markov's inequality,
\begin{align*}
& \P_0(\|\fmv-f_0\| > M_n\rho_n)\\
& \leq \P_0(\rho_n^2 < K^{-1}V_n) + \P_0(\|\fmv-f_0\|^2 \geq M_n^2 K^{-1}V_n) \\
& \lesssim \P_0(\rho_n^2 < K^{-1}V_n) + M_n^{-2} \frac{\P_0\|\fmv-f_0\|^2}{V_n}.
\end{align*}
This can be seen to vanish by applying Theorem~\ref{t:radius} to the term on the left in the upper bound, and by combining Lemma~\ref{l:loss} with the assumptions on $R_n$ for the term on the right. Statements \textit{1} and \textit{2} follow, the former by taking $M_n = 1$.

Since $\nu_n \leq 1$, the assumption $R_n \to \infty$ in statement \textit{3} is equivalent to $B_n/V_n\rightarrow\infty$.
Then letting $F : L^2(\calX,\mu) \to L^2(\calX,\mu)$ be the operator $f \mapsto \sum_{j=1}^m \nu_j \ip{f,\varphi_j} \varphi_j$, and id the identity operator, the assumption can be further re-written as
\begin{equation}\label{e:Rnzero}
V_n / \|(\id - F)f_0\|^2 \to 0.
\end{equation}
Note that by the triangle inequality,
\begin{align*}
& \P_0\big(\|\fmv - f_0\| \leq \rho_n\big) \\
& \leq \P_0\big(\|(\id - F) f_0\| \leq \|\fmv - Ff_0\| + M\rho_n\big) \\
& \leq \P_0\big(\|(\id - F)f_0\| \leq 2 \|\fmv - F f_0\|\big) + \P_0\big(\|(\id - F) f_0\| \leq 2M\rho_n\big).
\end{align*}
The second probability on the right hand side is seen to vanish as $n\to\infty$ by combining \eqref{e:Rnzero} with Theorem~\ref{t:radius}. Regarding the first probability, \eqref{e:pzerodelta} in the proof of Lemma~\ref{l:loss} below, gives
\[ \P_0\|\fmv - Ff_0\|^2 \lesssim n^{-1} \sum_{j=1}^m \nu_j^2 \]
so by Markov's inequality and \eqref{e:Rnzero} it follows that
\begin{align*}
\P_0\big(\|(\id - F)f_0\|
&\leq 2 \|\fmv - F f_0\|\big) \\
& \lesssim \frac{\P_0\|\fmv-Ff_0\|^2}{\|(\id - F)f_0\|^2} \\
& \lesssim \frac{n^{-1}\sum_{j=1}^m \nu_j^2}{V_n} \frac{V_n}{\|(\id - F)f_0\|^2} \to 0.
\end{align*}

\subsection{Proof of Lemma \ref{l:krr}}\label{sec:proof:krr}
Any $f \in \bb H_m$ is of the form $f(x) = \sum_{j=1}^m a_j h_j(x) = K_{x \bs u} \bs a$ for some $\bs a = (a_1,\ldots, a_m) \in \bb R^m$. Moreover, such an $f$ has squared RKHS norm
\[ \|f\|_{\bb H}^2 = \sum_j \sum_k a_j a_k \ip{h_j,h_k}_{\bb H} = \sum_j \sum_k a_j a_k \Pi u_j u_k = \bs a' K_\uu \bs a. \]
It follows that we need to minimise the objective function
\[ \bs a \mapsto (\bs y - K_\fu \bs a)'(\bs y - K_\fu \bs a) + \sigma^2 \bs a' K_\uu \bs a. \]
Setting the derivative equal to zero yields
\begin{equation}\label{e:score}
-K_\uf \bs y + (\sigma^2 K_\uu + K_\uf K_\fu)\bs a = 0.
\end{equation}
Since the second derivative $2(\sigma^2 K_\uu + K_\uf K_\fu)$ is positive definite, the solution $\bs a^\star = (\sigma^2 K_\uu + K_\uf K_\fu)^{-1} K_\uf \bs y$ minimises the objective function. It follows that
\[ \arg \min_{f \in \mathbb H_m} \sum_{i=1}^n (y_i-f(x_i))^2 + \sigma^2 \|f\|_{\bb H}^2 = K_{\cdot \bs u} \bs a^\star = \fmv. \qedhere \]

\subsection{Proof of Lemma \ref{l:loss}}\label{sec:proof:loss}
The proof follows similar lines of reasoning as \cite{bhattacharya:pati:yun:17} and \cite{hadji2022} in context of GP and distributed GP regression, using kernel ridge regression techniques. Here these standard techniques are adapted to the variational approximation.

Letting $\id$ denote the identity operator and $F$ denote the operator
\[ F  = \sum_{j=1}^m \nu_j \ip{\,\cdot\,,\varphi_j} \varphi_j \]
on $L^2(\calX,\mu)$, we obtain an identity for the sum of the bias terms
\[ \sum_{j=1}^m (1-\nu_j)^2\ip{f_0,\varphi_j}^2 + \sum_{j>m} \ip{f_0,\varphi_j}^2 = \|(\id-F)f_0\|^2. \]

Define $\Delta f_0 = \fmv - Ff_0$. Since $\|\fmv-f_0\|^2 \leq 2 \|\Delta f_0\|^2 + 2 \|(\id - F)f_0\|^2$, it suffices to show that
\begin{equation}\label{e:pzerodelta}
\P_0 \|\Delta f_0\|^2 \lesssim \frac1n\sum_{j=1}^m \nu_j^2.
\end{equation}
This is done by characterising the variational posterior mean $\fmv$ as the root of a ``score'' function. Let us write $\fmv = K_{x\bs u}\bs a^\star$ as in the proof of Lemma~\ref{l:krr}. Combining \eqref{e:score} with the identities \eqref{e:k}, it follows that $\fmv = \varphi_{1:m}(\cdot)'\Lambda \bs a^\star$ solves the equation
\begin{align}
\nonumber
\bs 0 & = \Lambda \Phi' \bs y - (\sigma^2 + \Lambda\Phi'\Phi)\Lambda \bs a^\star \\
\label{e:scoreu}
& = \Lambda\Phi'\left( \bs y - \bm{\fmv(x_1) \\ \vdots \\ \fmv(x_n)} \right) - \sigma^2 \bm{\ip{\fmv,\varphi_1} \\ \vdots \\ \ip{\fmv,\varphi_m}}
\end{align}
where we recall that $\Lambda = \operatorname{diag}(\lambda_1,\ldots,\lambda_m)$ and $\Phi_{i,j} = \varphi_j(x_i)$. Multiplying the $j$-th entry of the vector in \eqref{e:scoreu} with $\varphi_j$ gives that $\fmv$ is the root of the ``score'' operator $\hat S_n$
\[ \hat S_n (g) = \sum_{j=1}^m \Big[\frac 1n \sum_{i=1}^n (y_i - g(x_i)) \varphi_j(x_i)\Big]\lambda_j \varphi_j - \frac{\sigma^2}{n} g, \qquad g \in \bb H_m, \]
that is, $\hat S_n (\fmv) = 0$. Similarly, $Ff_0$ is the root of the ``population score''
\[ S_n (g) = \P_0 \hat S_n (g) = \sum_{j=1}^m \ip{f_0-g,\varphi_j} \lambda_j \varphi_j - \frac{\sigma^2}{n} g, \qquad g \in \bb H_m. \]
Let $f_{0,m} = \sum_{j=1}^m \ip{ f_0,\varphi_j}\varphi_j$ denote the orthogonal projection of $f_0$ onto $\bb H_m$. Then
\[ S_n (\fmv) = T f_{0,m} - \sum_{j=1}^m \frac{n\lambda_j + \sigma^2}{n} \ip{\fmv,\varphi_j} = T(f_{0,m} - F^{-1} \fmv), \]
so $- \Delta f_0 = Ff_0 - \fmv = FT^{-1}S_n(\fmv)$ and
\[ \|\Delta f_0\|^2 = \|FT^{-1} S_n(\fmv)\|^2 \leq 2 \|FT^{-1}(\hat S_n(Ff_0)+S_n(\fmv))\|^2 + 2 \|FT^{-1} \hat S_n(Ff_0)\|^2. \]
We establish at the end of this proof that
\begin{equation}\label{e:deltaleft}
\P_0 \|FT^{-1}(S_n(\fmv)+\hat S_n(Ff_0))\|^2 \lesssim m^4n^{-D}+ o(\P_0\|\Delta f_0\|^2),
\end{equation}
so the preceding two displays together yield
\begin{equation}\label{e:pdelta}
\P_0 \|\Delta f_0\|^2 \lesssim  \P_0 \|FT^{-1} \hat S_n(Ff_0)\|^2 + m^4n^{-D}+ o(\P_0\|\Delta f_0\|^2).
\end{equation}
Now we use that $y_i = f_0(x_i) + \varepsilon_i$ under $\P_0$, so
\begin{align*}
\hat S_n (Ff_0)
& = \hat S_n (Ff_0) - S_n (Ff_0) \\
& = \sum_{j=1}^m \Big[\frac1n \sum_{i=1}^n(y_i-Ff_0(x_i))\varphi_j(x_i) - \ip{(\id-F)f_0,\varphi_j}\Big]\lambda_j\varphi_j \\
&= \sum_{j=1}^m \Big[\frac1n \sum_{i=1}^n(f_0-Ff_0)(x_i)\varphi_j(x_i) - \ip{(\id-F)f_0,\varphi_j} + \frac1n\sum_{i=1}^n \varepsilon_i \varphi_j(x_i) \Big]\lambda_j\varphi_j
\end{align*}
so since the $\varepsilon_i$ are independent of the $x_i$ and both are i.i.d.,
\begin{align*} \P_0 \|FT^{-1}\hat S_n (Ff_0)\|^2 
& = \sum_{j=1}^m \nu_j^2 \P_0 \Big[\frac1n \sum_{i=1}^n[(\id-F)f_0](x_i)\varphi_j(x_i) - \ip{(\id-F)f_0,\varphi_j} + \frac1n\sum_{i=1}^n \varepsilon_i \varphi_j(x_i) \Big]^2 \\
& = \sum_{j=1}^m \nu_j^2 \frac1n \var_0\big([(\id-F)f_0](x_1)\varphi_j(x_1)\big)  + \frac{\sigma^2}{n}\sum_{j=1}^m \nu_j^2
 \lesssim \frac1n \sum_{j=1}^m \nu_j^2,
\end{align*}
the bottom line following from the boundedness of $f_0$. This together with \eqref{e:pdelta} for $D$ large enough yields \eqref{e:pzerodelta}, which concludes the proof.

~

\noindent\textsc{Proof of \eqref{e:deltaleft}.} Since $S_n(Ff_0)=\hat S_n(\fmv)=0$,
\begin{align*}
\hat S_n(Ff_0)+S_n(\fmv)  
& = \hat S_n(Ff_0) - \hat S_n(\fmv) + S_n(\fmv) - S_n(Ff_0) \\
& = \sum_{j=1}^m \Big[\frac 1n \sum_{i=1}^n (\Delta f_0)(x_i) \varphi_j(x_i) -  \ip{\Delta f_0,\varphi_j} \Big] \lambda_j \varphi_j
\end{align*}
and using $\Delta f_0 = \sum_{k=1}^m \ip{\Delta f_0,\varphi_k}\varphi_k$ and the Cauchy-Schwarz inequality, it follows that
\begin{align}
\nonumber
& \P_0 \| F T^{-1} ( \hat S_n (Ff_0)+S_n (\fmv) )\|^2 \\
\nonumber
& = \P_0 \Big\| \sum_{j=1}^m \Big[ \sum_{k=1}^m \ip{\Delta f_0,\varphi_k} \Big[ \frac 1n \sum_{i=1}^n  \varphi_k(x_i)\varphi_j(x_i) - \delta_{jk} \Big] \Big] \nu_j \varphi_j \Big\| \\
\nonumber
& = \sum_{j=1}^m \nu_j^2 \P_0 \Big[\sum_{k=1}^m \ip{\Delta f_0,\varphi_k} \Big[ \frac1n \sum_{i=1}^n \varphi_k(x_i) \varphi_j(x_i) - \delta_{jk}\Big]\Big]^2 \\
\label{e:deltaleftbound}
& \leq \P_0 \|\Delta f_0\|^2 \sum_{j=1}^m \nu_j^2 \sum_{k=1}^m \Big[ \frac1n \sum_{i=1}^n \varphi_k(x_i) \varphi_j(x_i) - \delta_{jk}\Big]^2.
\end{align}
Splitting over the events
\begin{equation}\label{e:Omega}
\Omega_{j,m} := \Big\{ \sup_{1\leq k\leq m} \Big[\frac1n\sum_{i=1}^n\varphi_j(x_i)\varphi_k(x_i)-\delta_{jk}\Big]^2 > C \frac{\log n}{n} \Big\}
\end{equation}
and their complements, and then using Assumption~\ref{a:bound}, the term in \eqref{e:deltaleftbound} is at most of the order
\begin{align*}
\frac{m \log n}{n} \sum_{j=1}^m \nu_j^2 \P_0\|\Delta f_0\|^2 + C_\varphi^2m^2 \max_{1 \leq j \leq m} \P_0 1_{\Omega_{j,m}} \|\Delta f_0\|^2.
\end{align*}
Now \eqref{e:deltaleft} follows from
\[ \frac{m \log n}{n} \sum_{j=1}^m \nu_j^2 \leq \frac{m^2 \log n}{n} \to 0 \]
(by assumption) and, as we prove now
\begin{equation}\label{e:f0negl}
\max_{1 \leq j \leq m} \P_01_{\Omega_{j,m}} \| \Delta f_0 \|^2 \lesssim m^2 n^{-D}.
\end{equation}
To this end, recall that from \eqref{e:vpmu}, $\fmv(x) = \varphi_{1:m}(x)'(\sigma^2\Lambda^{-1} + \Phi'\Phi)^{-1}\Phi'\bs y$, so
\begin{align*}
\|\fmv\|^2 = \sum_{k=1}^m \ip{\fmv,\varphi_k}^2 & = \operatorname{tr}\Big[(\sigma^2\Lambda^{-1}+\Phi'\Phi)^{-1} \Phi'\bs y \bs y' \Phi(\sigma^2\Lambda^{-1}+\Phi'\Phi)^{-1}\Big] \\
& \leq \lambda_1^2 \sigma^{-4} \bs y'\Phi\Phi'\bs y \lesssim \sum_{j=1}^m \Big[\sum_{i=1}^n y_i \varphi_j(x_i) \Big]^2 \lesssim m n \sum_{i=1}^n y_i^2,
\end{align*}
and hence, since $y_i = f_0(x_i)+\varepsilon_i$ for some bounded $f_0$ and $\varepsilon_i$ independent of $x_i$,
\[ \P_0 1_{\Omega_{j,m}} \|\fmv\|^2 \lesssim mn \P_01_{\Omega_{j,m}}\Big[\sum_{i=1}^n f_0(x_i)^2 + n\sigma^2\Big] \lesssim mn^2 \P_0(\Omega_{j,m}). \]
This together with $\|Ff_0\| \leq \|f_0\| < \infty$ yields
\[ \P_0 1_{\Omega_{j,m}} \| \Delta f_0 \|^2 \leq 2 \P_0(\Omega_{j,m}) \|Ff_0\|^2 + 2 \P_0 1_{\Omega_{j,m}} \|\fmv\|^2 \lesssim mn^2 \P_0(\Omega_{j,m}). \]
The inequality \eqref{e:f0negl} follows from the above by combining a union bound and Hoeffding's inequality:
\begin{equation*}
\P_0(\Omega_{j,m}) \leq m\cdot 2\exp\Big( \frac{-2n^2 C\log n}{4C_\varphi^2 n^2} \Big) = 2 m n^{-D-2},
\end{equation*}
with $D = C/2C_\varphi^2-2$, where the constant  $C$ can be chosen arbitrarily large.

\section{Proof of the corollaries}\label{sec:proof:cor}

\subsection{Proof of Corollary \ref{c:polycontr}}\label{sec:proof:polycontr}
Take $d/(d+2\alpha) \leq r < 1/2$. Note that $ m^2n^{-1} \log n=o(1)$, hence the conditions of Theorem~\ref{t:contraction} hold. Next note that
\begin{equation}\label{e:polyfd}
J_n = \max\{j : n\lambda_j \geq 1\} \asymp n^{d/(d+2\alpha)}\lesssim m,
\end{equation}
hence in view of \eqref{e:fdvar},
\begin{align}
V_n&\asymp  J_n/n+ \int_{m\wedge J_n}^\infty x^{-1-2\alpha/d} \,dx \asymp J_n/n+ J_n^{-2\alpha/d} \asymp n^{-2\alpha/(d+2\alpha)}\nonumber\\
W_n&\lesssim  J_n/n+ \int_{J_n}^{\infty} x^{-2-4\alpha/d} \,dx\asymp J_n/n+ n J_n^{-1-4\alpha/d} \asymp n^{-2\alpha/(d+2\alpha)}.\label{eq:var:poly}
\end{align}

To bound the bias term $B_n$, we consider \eqref{e:fdbias} and deal with the two terms on the right hand side separately. For the first term note
\begin{align}
\nonumber
\sum_{j=1}^{m \wedge J_n} (n\lambda_j)^{-2} \ip{f_0,\varphi_j}^2
& \lesssim \sum_{j=1}^{J_n} (nj^{-1-2\alpha/d})^{-2} \ip{f_0,\varphi_j}^2\\
\nonumber
& = n^{-2} \sum_{j=1}^{J_n} j^{2+4\alpha/d-2\beta/d} j^{2\beta/d} \ip{f_0,\varphi_j}^2 \\
\nonumber
& \lesssim n^{-2} \max_{1 \leq j < J_n} j^{2+4\alpha/d-2\beta/d} \|f_0\|_\beta^2 \\
\nonumber
& \lesssim (n^{-2} \vee n^{-2\beta/(d+2\alpha)}) \|f_0\|_\beta^2,
\end{align}
while for the second term we get
\begin{equation}\label{e:f0tail}
\sum_{j>m\wedge J_n}^\infty \ip{f_0,\varphi_j}^2 < (J_n\wedge m)^{-2\beta/d} \sum_{j>m\wedge J_n}^\infty j^{2\beta/d} \ip{f_0,\varphi_j}^2 \lesssim n^{-2\beta/(d+2\alpha)} \|f_0\|_\beta^2.
\end{equation}
hence the contraction rate in \eqref{eq:thm:contraction} is $\epsilon_n=n^{-(\beta\wedge\alpha)/(d+2\alpha)}$.

The sub-optimal contraction rate for insufficient amount of inducing variables is a direct consequence of Theorem \ref{t:badcontr}.
~

\subsection{Proof of Corollary \ref{c:polyuq}}\label{sec:proof:polyuq}
Note that $m^2n^{-1}\log n = n^{2r-1}\log n \to 0$, hence we can apply Theorem~\ref{t:coverage} to the present setting.

By \eqref{e:fdvar}, similarly to \eqref{eq:var:poly}, we get that
\begin{equation}\label{e:fdpolyvar}
V_n \asymp \frac{m \wedge n^{d/(d+2\alpha)}}{n} + m^{-2\alpha/d}\quad\text{and}\quad W_n\lesssim \frac{m \wedge n^{d/(d+2\alpha)}}{n}.
\end{equation}
Then similarly as in the proof of Corollary~\ref{c:polycontr} we get for $f_0 \in \calS^\beta$, that
\begin{align*}
\sum_{j=1}^m (1-\nu_j)^2\ip{f_0,\varphi_j}^2 & \lesssim (n^{-2} \vee n^{-2\beta/(d+2\alpha)}) \\
\sum_{j=m+1}^\infty \ip{f_0,\varphi_j}^2 & \lesssim m^{-2\beta/d} \|f_0\|_\beta^2.
\end{align*}
Therefore,
\begin{align}
R_n & = \frac{B_n+W_n}{V_n}
\lesssim \frac{n^{-1}(m \wedge n^{d/(d+2\alpha)}) + n^{-2\beta/(d+2\alpha)} + m^{-2\beta/d}}{n^{-1}(m \wedge n^{d/(d+2\alpha)})+m^{-2\alpha/d}}.\label{eq:UB:R_n}
\end{align}

In the first case, \eqref{eq:UB:R_n} is further bounded from above by a multiple of $m^{-2\beta/d}/m^{-2\alpha/d}=O(1)$. In the second case, \eqref{eq:UB:R_n} is $o(1)$. Finally, in the third case, let us take $f_0 \in \calS^\beta$ given by $f_0 = \sum_{j=1}^\infty j^{-1/2-q} \varphi_j$, for some $\beta/d<q<\alpha/d$. Then in view of \eqref{e:fdbias},
\[ \sum_{j=1}^{m \wedge J_n} (1-\nu_j)^2 \ip{f_0,\varphi_j}^2 + \sum_{j > m \wedge J_n} \ip{f_0,\varphi_j}^2 \gtrsim \sum_{j > m \wedge J_n} \ip{f_0,\varphi_j}^2 \asymp (m \wedge J_n)^{-2q} \]
hence by recalling \eqref{e:fdpolyvar},
\[ R_n \gtrsim \frac{m^{-2q} \vee n^{-2qd/(d+2\alpha)}}{m^{-2\alpha/d}+n^{-2\alpha/(d+2\alpha)}} \to \infty. \]
Therefore, in all three cases the statement follows from Theorem~\ref{t:coverage}.

Finally, the upper bound for the radius follows from Theorem \ref{t:radius} and \eqref{eq:var:poly} with $\alpha=\beta$.

~

\subsection{Proof of Corollary \ref{c:expcontr}}\label{sec:proof:expcontr}
The proof goes similarly to Corollary \ref{c:polycontr}. First assume that $m$ is lower bounded by
\begin{equation}\label{e:fdexp}
J_n := (\tau_n^{-1} \log n)^d = n^{d/(d+2\alpha)}
\end{equation}
and note that $\lambda_{n,J_n}  \asymp \exp(-\tau_n J_n^{1/d}) = n^{-1}$. Hence, similarly to \eqref{e:fdbias}
\begin{equation*}
B_n  \asymp \sum_{j=1}^{m\wedge J_n} (n\lambda_{n,j})^{-2}\ip{f_0,\varphi_j}^2 + \sum_{j>m \wedge J_n} \ip{f_0,\varphi_j}^2.
\end{equation*}
We deal with the two terms on the right hand side of the preceding display separately.

For any $f_0 \in \calS^\beta$, we have
\begin{align}
\sum_{j > m \wedge J_n} \ip{f_0,\varphi_j}^2 \leq J_n^{-2\beta/d} \sum_{j > m \wedge J_n} j^{2\beta/d} \ip{f_0,\varphi_j}^2 \lesssim n^{-2\beta/(d+2\alpha)} \|f_0\|_\beta^2  \label{eq:exp:B_n1}
\end{align}
and
\begin{align}
\sum_{j=1}^{m \wedge J_n}(n\lambda_{n,j})^{-2} \ip{f_0,\varphi_j}^2
& \leq \max_{1 \leq i \leq  J_n} (n\lambda_{n,i})^{-2}i^{-2\beta/d} \sum_{j=1}^\infty j^{2\beta/d} \ip{f_0,\varphi_j}^2\nonumber \\
& \lesssim n^{-2} \Big( \exp(2\tau_n) \vee \frac{\exp(2\tau_n J_n^{1/d})}{ J_n^{2\beta/d}}\Big) \|f_0\|_\beta^2 \nonumber\\
& \lesssim (n^{-2} \vee n^{-2\beta/(d+2\alpha)}) \|f_0\|_\beta^2 \label{eq:exp:B_n2}
\end{align}
where we used that the function $i\mapsto \exp(2\tau_n i)i^{-2\beta/d}$ is convex in $i$, so the maximum occurs at one of the endpoints.

Next we deal with the variance terms $W_n$ and $V_n$, defined in \eqref{def:bias:var} and \eqref{def:spread}, respectively. Similarly to \eqref{e:fdvar} and \eqref{eq:var:poly}, for $J_n$ given in \eqref{e:fdexp}, 
\begin{align}
V_n  \asymp \frac{m \wedge J_n}{n} + \sum_{j>m \wedge J_n} \lambda_{n,j}&\asymp \frac{(m \wedge J_n)}{n}+ \int_{m \wedge J_n}e^{-\tau_n x^{1/d}}dx\label{e:fdvarexp}\\
& \lesssim n^{-\frac{2\alpha}{d+2\alpha}} + \int_{(m \wedge J_n)^{1/d}} x^{d-1} e^{-\tau_n x} \,dx.\nonumber
\end{align}
By partial integration and induction we get that the rightmost integral satisfies
\begin{align}
 \int_{(m \wedge J_n)^{1/d}} x^{d-1} e^{-\tau_n x} \,dx\asymp  (\tau_n^{-1}(m \wedge J_n)^{1-1/d} \vee \tau_n^{-d}) \exp(-\tau_n(m \wedge J_n)^{1/d}).\label{e:expintbound}
\end{align}
Since $m \geq J_n = (\tau_n^{-1} \log n)^d$ it is further bounded by
\[ n^{d/(d+2\alpha)}(\log n)^{-1}\exp(-\log n) = n^{-2\alpha/(d+2\alpha)}/\log n. \]

Furthermore, by similar computations
\begin{align}
\nonumber
W_n &\asymp \frac{m \wedge J_n}{n} + \sum_{j>m \wedge J_n} \lambda_{n,j} \\
\nonumber
& \lesssim  n^{-2\alpha/(d+2\alpha)} + n\int_{(m \wedge J_n)^{1/d}} x^{d-1} e^{-2\tau_n x} \,dx\\
\nonumber
&\lesssim  n^{-2\alpha/(d+2\alpha)} + n (\tau_n^{-1}(m \wedge J_n)^{1-1/d} \vee \tau_n^{-d}) \exp(-2\tau_n(m \wedge J_n)^{1/d}) \\
&\lesssim n^{-2\alpha/(d+2\alpha)}.\label{eq:exp:W_n}
\end{align}
Our contraction rate result follows from Theorem~\ref{t:contraction} with $\epsilon_n^2=n^{-2(\beta\wedge\alpha)/(d+2\alpha)}$.

The sub-optimal contraction rate for insufficient amount of inducing variables is a direct consequence of Theorem \ref{t:badcontr}.

\subsection{Proof of Corollary \ref{c:expuq}}\label{sec:proof:expuq}

Note that by assumption we have $m^2n^{-1}\log n \to 0$, hence we can apply Theorem \ref{t:coverage}. Therefore it is sufficient to investigate the asymptotic behaviour of the fraction $R_n=(B_n+W_n)/V_n$, with the terms defined in \eqref{def:bias:var} and \eqref{def:spread}. We also recall the definition of $J_n$ given in \eqref{e:fdexp} for the exponentially decaying eigenvalues.

In the second case, following from the bounds \eqref{eq:exp:B_n1}, \eqref{eq:exp:B_n2} and \eqref{eq:exp:W_n}, the numerator of $R_n$ is $o(1)$. At the same time,  in view of assertions \eqref{e:fdvarexp}, \eqref{e:expintbound} and $ \tau_n m^{1/d} = n^{r/d-1/(d+2\alpha)} \log n \to 0$, the denominator is bounded from below as
\begin{align*}
V_n\gtrsim\sum_{j>m}\lambda_{n,j} \gtrsim \tau_n^{-d} \exp(-\tau_n m^{1/d}) \gtrsim \tau_n^{-d} \to \infty,
\end{align*}
and therefore $R_n=o(1)$.

Next, in the first case, we have $m \wedge J_n \asymp n^r \wedge J_n \asymp J_n$, hence again in view of  \eqref{eq:exp:B_n1}, \eqref{eq:exp:B_n2}, $B_n\lesssim n^{-2} \vee n^{-2\beta/(d+2\alpha)}$. Moreover, following from \eqref{e:fdvarexp} and \eqref{eq:exp:W_n}, we get that $V_n\gtrsim J_n/n = n^{-2\alpha/(d+2\alpha)}$ and $W_n\lesssim n^{-2\alpha/(d+2\alpha)}$, respectively. Combining this upper bounds results in that $R_n=O(1)$.

Finally, in the third case, as in the proof of Corollary~\ref{c:polyuq}, let $f_0$ be the function $\sum_{j=1}^\infty j^{-1/2-q} \varphi_j \in \calS^\beta$ for some $\beta /d < q < \alpha/ d$. In view of \eqref{e:fdbias}
\[B_n \gtrsim \sum_{j > m \wedge J_n} \ip{f_0,\varphi_j}^2 \asymp (m \wedge J_n)^{-2q} \gtrsim n^{-2qd/(d+2\alpha)}. \]
Furthermore, \eqref{e:fdvarexp} and \eqref{e:expintbound} together with  $m \geq J_n$ imply that,
\[ V_n \lesssim J_n/n + \tau_n^{-1}J_n^{1-1/d} \exp(-\tau_n J_n^{1/d}) \lesssim n^{-2\alpha/(d+2\alpha)}. \]
Therefore, $R_n\gtrsim n^{-2qd/(d+2\alpha)}/n^{-2\alpha/(d+2\alpha)}\to \infty$. We conclude that, in all three cases the statements now follow directly from Theorem~\ref{t:coverage}.

Finally, the upper bound for the radius follows from Theorem \ref{t:radius} and the upper bound for $V_n$ in the proof of Corollary \ref{sec:proof:expcontr} with $\alpha=\beta$.

\section{Variational posterior for general inducing variables}\label{sec:vb:ind:var}
Recall that the posterior is approximated by the variational distribution
\begin{equation}\label{e:inducedPosterior}
\Psi = \int \Pi(\cdot | \bs u) \,d\Psi_{\bs u}(\bs u)
\end{equation}
on $L^2(\calX,\mu)$. The above display is equivalent to
\begin{equation}\label{e:VP}
\frac{d\Psi}{d\Pi}(f) = \frac{d\Psi_{\bs u}}{d\Pi_{\bs u}}(\bs u(f)),
\end{equation}
where $\Pi_{\bs u}$ is the distribution of the vector of inducing variables ${\bs u}$ under the prior $\Pi$ (this follows from existence of the Radon-Nikodym derivative on the right, which in turn is guaranteed by the assumption that $\Psi_{\bs u}$ is a non-degenerate Gaussian). The identity \eqref{e:VP} shows that the variational posterior $\Psi$ is parameterised by $f$ through ${\bs u}$ only.

Then, in view of Bayes' theorem
\[ \frac{d\Pi(\,\cdot\,|\bs x,\bs y)}{d\Pi} = \frac{p_f(\bs x,\bs y)}{p(\bs x,\bs y)} \]
and some additional elementary algebraic manipulations
\begin{align*}
& \KL(\Psi \| \tp) \\
& = \int \log \frac{d\Psi}{d\tp} \,d\Psi \\
& = \int \log \frac{d\Psi}{d\Pi} \,d\Psi + \int \log \frac{d\Pi}{d\tp} \,d\Psi \\
& = \int \log \frac{d\Psi_{\bs u}}{d\Pi_{\bs u}} \,d\Psi_{\bs u} - \iint \log p_f(\bs x, \bs y) \,d\Psi(f|{\bs u})\,d\Psi_{\bs u}({\bs u}) + \log p(\bs x, \bs y).
\end{align*}
In the variational procedure, we aim to minimise this with respect to $\Psi$. Using that the KL-divergence is non-negative, and $d\Psi(f|{\bs u}) = d\Pi(f|{\bs u})$ by construction, we obtain
\[ \log p(\bs x, \bs y) \geq \iint \log p_f(\bs x, \bs y) \,d\Pi(f|{\bs u})\,d\Psi_{\bs u}({\bs u})  - \int \log \frac{d\Psi_{\bs u}}{d\Pi_{\bs u}} \,d\Psi_{\bs u}. \]
The right-hand side is commonly referred to as evidence lower bound (ELBO). Minimising the KL-divergence is equivalent to maximising the ELBO. Note that by Jensen's inequality,
\begin{align}
\label{e:ELBO0}
\mathrm{ELBO}
& = \int \log \Big[ \frac{d\Pi_{\bs u}}{d\Psi_{\bs u}}({\bs u}) \, {\textstyle \exp\Big(\int \log p_f(\bs x, \bs y) \,d\Pi(f|{\bs u})\Big)} \Big] \,d\Psi_{\bs u}({\bs u}) \\
\label{e:ELBO}
& \leq \log \int {\textstyle \exp\Big(\int \log p_f(\bs x, \bs y) \,d\Pi(f|{\bs u})\Big)} \,d\Pi_{\bs u}({\bs u}),
\end{align}
and the maximum is attained by the distribution $\Psi_{\bs u}^*$ defined through
\begin{equation}\label{e:optimalVP}
\frac{d\Psi_{\bs u}^*}{d\Pi_{\bs u}}({\bs v}) = \frac{\exp\Big(\int \log p_f(\bs x, \bs y) \,d\Pi(f|{\bs u}={\bs v})\Big)}{\int \exp\Big(\int \log p_f(\bs x, \bs y) \,d\Pi(f|{\bs u})\Big)\,d\Pi_{\bs u}(\bs u)}.
\end{equation}

In general there is no guarantee that the maximizer  $\Psi_{\bs u}^*$ of the ELBO is tractable, but in case of the considered Gaussian process regression model and Gaussian variational class it has an explicit form. In view of \eqref{e:optimalVP},
\begin{align}
\nonumber
\frac{d\Psi_{\bs u}^*}{d\Pi_{\bs u}}({\bs u})
& \propto \exp\Big(-\frac{1}{2\sigma^2}\int\sum_{i=1}^n (y_i-f(x_i))^2 \,d\Pi(f|{\bs u})\Big) \\
\nonumber
& = \exp\Big(-\frac{1}{2\sigma^2}\sum_{i=1}^n \big[(y_i - \Pi(f(x_i)|{\bs u}))^2 + \operatorname{var}_{\Pi(\cdot|{\bs u})}(f(x_i))\big]\Big) \\
\nonumber
& \propto \exp\Big(-\frac{1}{2\sigma^2}\sum_{i=1}^n(y_i-\Pi(f(x_i)|{\bs u}))^2 \Big) \\
\label{e:uVarDist}
& = \exp\Big(-\frac{1}{2\sigma^2}\sum_{i=1}^n(y_i-K_{x_i \bs u} K_{{\bs u}{\bs u}}^{-1} \bs u)^2 \Big), 
\end{align}
where $K_{x{\bs u}} = \operatorname{cov}_{\Pi}(f(x),{\bs u}) = \Pi f(x){\bs u}^{\text T}$ and $K_{{\bs u}{\bs u}} = \operatorname{cov}_\Pi({\bs u},{\bs u}) = \Pi {\bs u}{\bs u}^{\text T}$, and we used that
\begin{equation*}\label{e:varGivenu}
\operatorname{var}_{\Pi(\cdot|{\bs u})}(f(x)) = k(x,x) - K_{x{\bs u}}K_{{\bs u}{\bs u}}^{-1}K_{{\bs u}x},
\end{equation*}
which is constant as a function of ${\bs u}$. 

Completing the square in \eqref{e:uVarDist}, it follows that the optimal variational distribution $\Psi_{\bs u}^*$ of ${\bs u}$ is Gaussian with mean $K_{{\bs u}{\bs u}}(\sigma^2 K_{{\bs u}{\bs u}} + K_{{\bs u}\bs f}K_{\bs f {\bs u}})^{-1} K_{{\bs u}\bs f}\bs y$ and covariance matrix $K_{{\bs u}{\bs u}}(K_{{\bs u}{\bs u}} + \sigma^{-2} K_{{\bs u}\bs f}K_{\bs f {\bs u}})^{-1} K_{{\bs u}{\bs u}}$. By \eqref{e:inducedPosterior}, the variational distribution $\Psi^*$ of $f$ is Gaussian, with mean function
\begin{align*}
\fmv(x) & = \int \Pi(f(x)|{\bs u}) \,d\Psi_{\bs u}^*({\bs u}) \\
& = K_{x{\bs u}} K_{{\bs u}{\bs u}}^{-1} \int {\bs u} \,d\Psi_{\bs u}^*({\bs u})\\
& = K_{x{\bs u}} (\sigma^2 K_{{\bs u}{\bs u}} + K_{{\bs u}\bs f}K_{\bs f {\bs u}})^{-1} K_{{\bs u} \bs f} \bs y.
\end{align*}
Furthermore, the covariance function is
\begin{align*}
(x,y) & \mapsto \int \Pi((f-\fmv)(x)(f-\fmv)(y)|{\bs u}) \,d\Psi^*_{\bs u}({\bs u}) \\
& = \int \cov_{\Pi(\cdot|{\bs u})}(f(x),f(y)) \,d\Psi^*_{\bs u}({\bs u}) + \cov_{\Psi^*_{\bs u}} (\Pi(f(x)|{\bs u}),\Pi(f(y)|{\bs u})) \\
& = k(x,y) - K_{x{\bs u}}K_{{\bs u}{\bs u}}^{-1} K_{{\bs u}y} + K_{x{\bs u}} K_{{\bs u}{\bs u}}^{-1} \cov_{\Psi_{\bs u}^*}({\bs u},{\bs u}) K_{{\bs u}{\bs u}}^{-1} K_{{\bs {\bs u}}y} \\
& = k(x,y) - K_{x{\bs u}}K_{{\bs u}{\bs u}}^{-1} K_{{\bs u}y} + K_{x{\bs u}} (K_{{\bs u}{\bs u}} + \sigma^{-2} K_{{\bs u}\bs f}K_{\bs f {\bs u}})^{-1} K_{{\bs u}y}.
\end{align*}

\newpage
\bibliographystyle{imsart-number}
\bibliography{references}
\end{document}